\documentclass[article,onefignum,onetabnum]{siamart171218}
\usepackage{amsfonts}
\usepackage{graphicx}
\usepackage{placeins}
\usepackage{comment}
\usepackage{MnSymbol}
\usepackage{enumerate}

\newcommand{\bn}{\mathbf{n}}

\newcommand{\e}{\varepsilon}
\newcommand{\ov}{\overline}
\newcommand{\R}{\mathbb R}

\newcommand{\rr}{{\bf r}}

\newcommand{\N}{{\bf N}}

\newcommand{\beq}{\begin{equation}}
\newcommand{\eeq}{\end{equation}}

\newsiamremark{remark}{Remark}
\newsiamremark{hypothesis}{Hypothesis}
\crefname{hypothesis}{Hypothesis}{Hypotheses}
\newsiamthm{claim}{Claim}

\headers{Well-posedness of  Cahn-Hilliard model for surface diffusion}{Niu, Xiang,Yan}

\title{Well-posedness of a modified degenerate Cahn-Hilliard model for surface diffusion}

\author{Xiaohua Niu\thanks{School of Mathematics and Statistics, Xiamen University of Technology, Xiamen, China
  (\email{xniu@connect.ust.hk}).}
\and Yang Xiang \thanks{Department of Mathematics, The  Hong Kong University of Science and Technollogy, Clear Water Bay, Kowloon, Hong Kong
(\email{maxiang@ust.hk}).}
\and Xiaodong Yan\thanks{Department of Mathematics, The University of Connecticut, Storrs, CT, USA
  (\email{xiaodong.yan@uconn.edu}).}
}

\begin{document}

\maketitle

\begin{abstract}
We study the well-posedness of a modified  degenerate Cahn-Hilliard type model for surface diffusion.
With degenerate phase-dependent diffusion mobility and additional stabilizing function, this model is
able to give the correct sharp interface limit.  We introduce a notion of weak solutions for the nonlinear model. The existence result is obtained by approximations  of the proposed model with nondegenerate mobilities.  We  also employ this method to prove existence of weak solutions to a related model where the chemical potential contains a nonlocal term originated from self-climb of dislocations in crystalline materials.
\end{abstract}



\section{Introduction}\label{sec:intro}

We consider the following modified degenerate Cahn-Hilliard type model
\begin{eqnarray}
\label {eqn:ch} g(u)\partial_t u&=&\nabla \cdot (M(u)\nabla \frac{\mu}{g(u)}), \ \   x \in \Omega \subset \mathbb{R}^2, t\in [0,\infty) \\
\label{eqn:chem} \mu&=&-\Delta u+\frac{1}{\e^2}q'(u).
\end{eqnarray}

When $g\equiv 1$, \eqref{eqn:ch}-\eqref{eqn:chem} becomes Cahn-Hilliard (CH) equation   with degenerate mobility. The degenerate Cahn-Hilliard equation has been widely studied  as a diffuse-interface model for phase separation in binary system \cite{CahHil58,CahHil71,CahTay94,DaiDu12, DaiDu16,EllGar96}. Over the years, the interface motion in the sharp limit has caught a lot of attention for various choice of mobility $M(u)$ and homogeneous free energy $q(u)$.
When $M(u) = 1-u^2$ and $q$ being
either the logarithmic free energy $$q(u)=\frac{\theta}{2}\left[(1+u)\ln(1+u)+(1-u)\ln(1-u)\right]+\frac{1}{2}(1-u^2)$$ with temperature $\theta=O(\e^{\alpha})$ or
 the double obstacle potential
$$q(u) = 1-u^2 \text{ for  } |u|\leq 1, \text {    } q(u)=\infty \text { otherwise},$$
Cahn, Elliott, and Novick-Cohen [19] showed via asymptotic expansions that the
sharp-interface limit in the time scale $O(\e^{-2})$ is interface motion by surface diffusion. Sharp interface limits for different time scales were discussed in \cite{DaiDu12} for highly disparate diffusion mobility $M(u)=1+u$ and smooth double well $q(u)=\frac{1}{4}(1-u^2)^2$. In particular, the system evolves in $t=O(\e^{-2})$ time scale according to the combination of a one-sided modified Mullins–Sekerka problem in the phase with nonzero constant mobility and a
nonlinear diffusion process that solves a quasi-stationary porous medium equation in the phase with
small mobility. A later work by the same authors \cite{DaiDu14} derived sharp interface limit for $O(\e^{-2})$ time scale with  diffusion mobility $M(u)=|1-u^2|$ and smooth double well potential $q(u)$,  noting the effect of the diffusion field on the interface motion  as a jump of fluxes. The analysis was done on the
(unphysical) solution branch with $|u|>1$ on some region. For $M(u)=1-u^2$ and $q=\frac{1}{4}(1-u^2)^2$, Lee, M\"unch and S\"uli  \cite{LeeMunSul15} considered the physical branch of solution
where $|u|<1$ everywhere  and showed that there is an additional nonlinear bulk diffusion term appearing to leading order of the sharp interface limit.  Further study in \cite{LeeMunSul16} indicates that  the
leading order sharp-interface motion depends sensitively on the choice of mobility.

The existence of weak solutions for degenerate Cahn-Hilliard equation was proved by Elliotte and Garcke \cite{EllGar96} (see \cite{Yin92} for 1D case). Their results include the case $M(u) = 1-u^2$ and $q$ being
the logarithmic free energy. Dai and Du \cite{DaiDu16} introduced a different notion of weak solutions for degenerate Cahn-Hilliard equation with mobility $M(u)=|1-u^2|^m$ and smooth double well potentials; they showed that their model accommodates the  Gibbs-Thomson effect, which was not by the method in \cite{EllGar96}.

There is a critical issue in modeling surface diffusion by the degenerate Cahn-Hilliard model \cite{GugSpaKas08,RRV06}, due to the presence of incompatibility in the asymptotic matching between the outer and inner expansions. R\"{a}tz, Ribalta, and Voigt
(RRV) \cite{RRV06} fixed this incompatibility by introducing a singular factor $1/g\left( u\right) $ in front of the chemical potential $\mu$ to force it to vanish in the far field. Their model essentially consists of equations \eqref{eqn:ch}-\eqref{eqn:chem} without the $g(u)$ term on the left side of \eqref{eqn:ch}, and other terms for modeling   heteroepitaxial growth
of thin films. The RRV model with the stabilizing function $g(u)$ has been validated by numerical simulations  \cite{RRV06} and asymptotic analyses \cite{GugSpaKas08,RRV06}. It has been successfully generalized to many applications, e.g.,  growth of nanoscale membranes \cite{Voigt2018}, dewetting of ultrathin films \cite{Naffouti2017}, and grain boundary formation in nanoporous metals \cite{Voorhees2021}. Recently, a phase field model for dislocation self-climb by vacancy pipe diffusion was developed based on degenerate Cahn-Hilliard model with such stabilizing function \cite{NXY20}.
However, to the best of our knowledge,  well-posedness of
these degenerate  Cahn-Hilliard models with singular factor that give the correct sharp interface limit for surface diffusion has not been established in the literature.


In this paper, in order to prove the well-posedness of
the RRV type Cahn-Hilliard model with correct sharp interface limit for surface diffusion, we propose a modified degenerated Cahn-Hilliard model as given in \eqref{eqn:ch}-\eqref{eqn:chem}, and discuss its well-posedness  and sharp interface limit. In particular, we have modified the original RRV model so that the equation can be written in the form of gradient flow of the total energy.

Our first result is  a sharp interface limit equation for \eqref{eqn:ch} and \eqref{eqn:chem} via formal asymptotic analysis. We obtain the following sharp interface equation
\begin{equation}
v=\lambda \partial_{ss}(\alpha \kappa)
\end{equation}
as $\e \rightarrow \infty$.  Here $\lambda<0$, $\alpha<0$  are constants whose exact forms are derived in  section \ref{sec:asym}. This validates this equation as a diffuse-interface model for surface diffusion.

Our main  result concerns the well-posedness of the initial value problem of  \eqref{eqn:ch}-\eqref{eqn:chem}. For this purpose,
 we set $\Omega=[0,2\pi]^n$ and consider the following problem  in a periodic setting when  $n\leq 2$.
\begin{eqnarray}
\label{eqn:u} g(u)\partial_t u&=&\nabla \cdot (M(u)\nabla \frac{\mu}{g(u)}), \hspace{0.4in} \text{ for } x \in \Omega, t\in [0,\infty) \\
\label{eqn:mu} \mu&=&-\Delta u+q'(u).
\end{eqnarray}
Here $g(u)=|1-u^2|^m$ for $2\leq m<\infty $,  $M(u)=M_0g(u)$ for some constant $M_0>0$ and $q(u)$ satisfies the following assumptions.
\begin{itemize}
 \item[(i)] $q(u)\in C^2(\R,\R)$ and there exist constants $C_i>0$, $i=1,\cdots,10$ such that for all $u \in \R$ and some $1\leq r <\infty$, the following growth assumptions hold.
\begin{eqnarray}\label{eqn:qu}
C_1|u|^{r+1}-C_2\leq q(u) &\leq &C_3|u|^{r+1}+C_4, \\ \label{eqn:q'u}|q'(u)|
&\leq &C_5|u|^r+C_6,\\ \label{eqn:q''u}C_t|u|^{r-1}-C_8&\leq& q''(u) \leq C_9|u|^{r-1}+C_{10}.
\end{eqnarray}

\end{itemize}
We see that the classical double well potential $q(u)=(1-u^2)^2$ satisfies \eqref{eqn:qu}-\eqref{eqn:q''u} with $r=3$.

Our existence proof is obtained via approximations  of the proposed model \eqref{eqn:u}-\eqref{eqn:mu} with positive mobilities. Given any $\theta >0$, we define

\begin{equation}
\label{eqn:Mtheta} M_{\theta}(u):= M_0 g_{\theta}(u)
\end{equation}
with
\begin{equation}
\label{eqn:gtheta} g_{\theta}(u):= \left\{ \begin{array}{cl}
                                 |1-u^2|^m & \text{  if  }|1-u^2|> \theta,\\
                                  \theta^m &\text { if } |1-u^2|\leq \theta.
                                 \end{array} \right.
\end{equation}

Our first step is to find a sufficiently regular solution for \eqref{eqn:u}-\eqref{eqn:mu} with mobility $M_{\theta}(u)$ and stablizing function $g_{\theta}(u)$ together with a smooth potential $q(u)$.
\begin{theorem}\label{thm:ndg}
Let  $M_{\theta},  g_{\theta}$ be defined by \eqref{eqn:Mtheta} and \eqref{eqn:gtheta}, under the assumptions \eqref{eqn:qu}-\eqref{eqn:q''u}, for any  $u_0 \in H^1(\Omega)$ and  any $T>0$, there exists a function $u_{\theta}$ such that
\begin{itemize}
\item[a)] $u_{\theta}\in L^{\infty}(0,T;H^1(\Omega))\cap C([0,T];L^p(\Omega))\cap L^2(0,T;W^{3,s}(\Omega))$, where $1\leq p <\infty$, $1\leq s <2$,
\item[b)]$\partial_tu_{\theta}\in L^2(0,T;(W^{1,q}(\Omega))')$ for $q>2$,
\item[c)]$u_{\theta}(x,0)=u_0(x) $ for all $x \in \Omega$,
\end{itemize}
which satisfies \eqref{eqn:u}-\eqref{eqn:mu} in the following weak sense
\begin{eqnarray}\notag
&&\int_0^T<\partial_tu_{\theta},\phi>_{(W^{1,q}(\Omega))',W^{1,q}(\Omega)} dt\\
&&=-\int_0^T\int_{\Omega}M_{\theta}(u_{\theta})\nabla \frac{-\Delta u_{\theta}+q'(u_{\theta})}{g_{\theta}(u_{\theta})}\cdot \nabla \frac{\phi}{g_{\theta}(u_{\theta})}dxdt
\label{eqn:thetandg}
\end{eqnarray}
for all $\phi\in L^2(0,T;W^{1,q}(\Omega))$ with $q>2$. In addition, the following energy inequality holds for all $t>0$.
\begin{eqnarray}
\notag
&&\int_{\Omega}\left (\frac{1}{2}|\nabla u_{\theta}(x,t)|^2+q(u_{\theta}(x,t))\right) dx\\
\notag &&+\int_0^t\int_{\Omega}M_{\theta}(u_{\theta}(x,\tau)\left|\nabla \frac{-\Delta u_{\theta}(x,\tau)+q'(u_{\theta}(x,\tau)}{g_{\theta}(u_{\theta}(x,\tau))}\right|^2dxd\tau \\
\label{eqn:engineq}&\leq&\int_{\Omega}\left(\frac{1}{2}|\nabla u_0(x)|^2+q(u_{0}(x))\right) dx.
\end{eqnarray}
\end{theorem}

Proof of theorem \ref{thm:ndg} is based on Galerkin approximations. Due to the presence of the stablizing function $g_{\theta}$, it is not obvious how to pass to the limit in the nonlinear term of the Galerkin approximations. Our key observation in this step is  strong convergence of  $\nabla u^N$ (up to a subsequence) in $L^2(\Omega_T)$ where $\Omega_T=\Omega\times [0,T])$ which allows us to pass to the limit in the nonlinear term.

To obtain the weak solution to \eqref{eqn:u}, we consider the limit of $u_{\theta_i}$ for a sequence $\theta_i\downarrow 0$.   The key challenge is how to pass to the limit in both sides of \eqref{eqn:thetandg}. In the degenerate Cahn-Hilliard case, the estimates for the positive mobility approximations yield a uniform bound for $\partial_t u_{\theta_i}$ and it is straighforward to pass to the limit on the left hand side in the approximating equations.  Moreover, the bound on $\partial_t u_{\theta_i}$, together with bound on $u_{\theta_i}$ yields strong  convergence of $\sqrt{M_i(u_{\theta_i})}$ in $C(0,T;L^n(\Omega))$. By this and the weak convergence of $\sqrt{M_i(u_{\theta_i})}\nabla \mu_{\theta_i}$ in $L^2(\Omega_T)$, Dai and Du \cite{DaiDu16} showed (up to a subsequence) that  $M_{\theta_i}(u_{\theta_i})\nabla \mu_{\theta_i}\rightharpoonup \sqrt{M(u)}\xi $ weakly in $L^2(0,T;L^{\frac{2n}{n+2}}(\Omega))$ where $\xi$ is the weak limit of $\sqrt{M_i(u_{\theta_i})}\nabla \mu_{\theta_i}$.  The main task  left  is  to show $\sqrt{M(u)}\xi =M(u)(-\nabla \Delta u+q''(u)\nabla u)$ and the limit equation becomes a weak form Cahn-Hilliard equation.  They  \cite{DaiDu16}  proved that this is almost true in the set where $u\neq \pm 1$. Their main idea is the following.  For small numbers $\delta_j$ monotonically decreasing to $0$, they  consider the limit in  a subset  $B_j$ of $\Omega_T$ where   approximate solutions converges uniformly and $|\Omega_T\backslash B_j|<\delta_j$.  By decomposing $B_j=D_j\cup\tilde D_j$ where mobility is bounded from below uniformly in $D_j$ and controlled above in $\tilde D_j$ by suitable multiples of $\delta_j$, they obtain the weak form equation for the limit function by passing  to the limit of $u_{\theta_i}$ on $D_j$ then letting $j$ goes to $\infty$.  Under further regularity assumptions on $\nabla \Delta u$, they obtained the explicit expression for $\xi$  in the weak form of the equation.

In this paper, we adapt their idea to our model.   There are two main difficulties. The first obtacle is the bound estimate on $\partial_t u_{\theta_i}$ blows up when $\theta_i$ goes to zero and we can not pass to the limit on the left hand side of \eqref{eqn:thetandg}; secondly, due to the presence of the stablizing function $g$ on the right hand side, it is more complicated to  derive an explicit expression of the weak limit of ${M_i(u_{\theta_i})}\nabla \frac{\mu_{\theta_i}}{g_{\theta_i}(u_{\theta_i})}$ in terms of $u$ on the right hand side  of the limit equation. To overcome the first difficulty, we derive an alternative form of \eqref{eqn:thetandg} by multiplying $g_{\theta}(u_{\theta})$ to both sides (valid due to regularity of $u_{\theta}$, c.f. section \ref{sec:regularity-u theta} and equation \eqref{eqn:utheta2}). From this, we obtain uniform estimates on $g_{\theta_i}(u_{\theta_i})\partial_tu_{\theta_i}$ which enables  us to pass to the limit on the left hand side of the alternate equation \eqref{eqn:utheta2}. To find limit form on the right hand side of \eqref{eqn:utheta2}, we need convergence of $\sqrt{M_{\theta_i}(u_{\theta_i})}$ in $C(0,T;L^p(\Omega))$. Due to the lack of control on $\partial_tu_{\theta_i}$, such convergence can not be derived directly using  Aubin-Lions Lemma \cite{DaiDu16}. Instead, we apply Aubin-Lions lemma to $G_i(u)=\int_0^u g_{\theta_i}(s) ds$ and derive convergence of $g_{\theta_i}(u_{\theta_i})$ (consequently on $M_{\theta_i}(u_{\theta_i})$) from convergence of $G_i$ through characterization of compact sets \cite{Sim86} in $L^p[0,T;B]$. We then follow the idea in \cite{DaiDu16} to pass to the limit on the right hand side of \eqref{eqn:utheta2}. Finally, we  identify an explicit expression of the weak limit of $\nabla \frac{\mu_{\theta_i}}{g_{\theta_i}(u_{\theta_i})}$ in terms of the weak limit $u$ under additional integrability assumptions on derivatives of $u$.

\begin{theorem}\label{thm:dg}
For any $u_0\in H^1(\Omega)$ and $T>0$, there exists a function $u:\Omega_T=\Omega\times [0,T]\rightarrow \R$ satisfying
\begin{itemize}
\item[i)] $u\in L^{\infty}(0,T;H^1(\Omega))\cap C([0,T];L^s(\Omega))$, where $1\leq s <\infty$,
\item[]
\item[ii)] $g(u)\partial_tu\in L^{p}(0,T;(W^{1,q}(\Omega))')$ for $1\leq p<2$ and $q>2$.
\item[]
\item[iii)] $u(x,0)=u_0(x)$ for all $x\in \Omega$,
\end{itemize}
which solves \eqref{eqn:u}-\eqref{eqn:mu} in the following weak sense
\begin{itemize}
\item[a)]There exists a set $B\subset \Omega_T$ with $|\Omega_T\backslash B|=0$ and a function $\zeta:\Omega_T\rightarrow \R^n$ satisfying $\chi_{B\cap P}M(u)\zeta\in L^{\frac{p}{p-1}}(0,T;L^{\frac{q}{q-1}}(\Omega,\R^n))$ such that
          \begin{equation}
           \int_0^T<g(u)\partial_tu,\phi>_{((W^{1,q}(\Omega))',W^{1,q}(\Omega))}dt=-\int_{B\cap P}M(u)\zeta \cdot \nabla \phi dxdt \label{eqn:dg}
          \end{equation} for all $\phi \in L^p(0,T;W^{1,q}(\Omega))$ with $p,q>2$. Here $P:=\{(x,t)\in \Omega_T: |1-u^2|\neq 0\}$ is the set where $M(u), g(u)$ are nondegenerate and $\chi_{B\cap P}$ is the characteristic function of set $B\cap P$.

\item[b)] Assume $u\in L^2(0,T;H^2(\Omega)).$  For any open set $U\in \Omega_T$ on  which $g(u)>0 $ and $\nabla \Delta u \in L^p(U)$ for some $p>1$, we have
\begin{equation}
\zeta=\frac{-\nabla \Delta u+q''(u)\nabla u}{g(u)}-\frac{g'(u)}{g^2(u)}\left(-\Delta u+q'(u) \right)\nabla u  \text{ a.e. in } U. \label{eqn:zetaU}
\end{equation}
\end{itemize}
Moreover, the following energy inequality holds for all $t>0$
\begin{eqnarray}
\notag &&\int_{\Omega}\left(\frac{1}{2}|\nabla u(x,t)|^2+q(u(x,t))\right)dz+\int_{\Omega_r\cap B\cap P}M(u(x,\tau))|\zeta(x,\tau)|^2 dxd\tau\\
\label{eqn:ineq}&\leq &\int_{\Omega}\left(\frac{1}{2}|\nabla u_0(x)|^2+q(u_0(x))\right)dx.
\end{eqnarray}
\end{theorem}

Adding an additional term $\frac{1}{\e}(-\Delta)^{\frac{1}{2}}u$ to the chemical potential in \eqref{eqn:chem}, we can apply the same method to derive existence of weak solutions of the modified model (see section \ref{sec:pf} for further details). Such nonlocal model originates from the phase field model for self-climb of dislocation loops \cite{NXY20}.

The paper is organized as follows. We shall derive sharp interface limit for \eqref{eqn:ch} and \eqref{eqn:chem} through formal asymptotic expansions in section \ref{sec:asym}. Section \ref{sec:ndg} is devoted to the proof of Theorem \ref{thm:ndg} and Theorem\ref{thm:dg}  is proved in section \ref{sec:dg}. Similar existence theorems for the modified model with an additional nonlocal term added to the chemical potential are presented in section \ref{sec:pf}.

\section{Sharp interface limit via asymptotic expansions }
\label{sec:asym}
In this section, we perform a formal asymptotic analysis to obtain the sharp interface limit of the proposed phase field model \eqref{eqn:ch}- \eqref{eqn:chem}as $\e \rightarrow 0$.
\subsection{Outer expansions}
We first perform expansion in the region far from the dislocations. Assume the expansion for $u$ is
\begin{equation}
u(x,y,t)=u^{(0)}(x,y,t)+u^{(1)}(x,y,t)\e+u^{(2)}(x,y,t)\e^2+\hdots
\end{equation}
Correspondingly, we have
\begin{eqnarray*}
&&M(u)=M(u^{(0)})+M'(u^{(0)})u^{(1)}\e +\left(M'(u^{(0)})u^{(2)}+\frac{1}{2}M''\left (u^{(0)}\right)\left(u^{(1)}\right)^2\right)\e^2+\hdots,\\
&&g(u)=g(u^{(0)})+g'(u^{(0)})u^{(1)}\e +\left(g'(u^{(0)})u^{(2)}+\frac{1}{2}g''(u^{(0)})\left(u^{(1)}\right)^2\right)\e^2+\hdots,\\
&&q'(u)=q'(u^{(0)})+q''(u^{(0)})u^{(1)}\e +\left(q''(u^{(0)})u^{(2)}+\frac{1}{2}q^{(3)}(u^{(0)})\left(u^{(1)}\right)^2\right)\e^2+\hdots,
\end{eqnarray*}
We also expand chemical potential $\mu$ as
\begin{equation}
\mu=\frac{1}{\e^2}\left(\mu^{(0)}+\mu^{(1)}\e+\mu^{(2)}\e^2+\hdots\right). \label{eqn:muexp}
\end{equation}
and rewrite equation \eqref{eqn:ch} as
\begin{equation}
g(u)\partial_t u=M_0\nabla \cdot (\nabla \mu-\mu \frac{g'(u)}{g(u)}\nabla u). \label{eqn:ch1}
\end{equation}
Set
$$
w=-\mu \frac{g'(u)}{g(u)}=\frac{1}{\e^2}\left(w^{(0)}+w^{(1)}\e+w^{(2)}\e^2+\hdots\right).
$$
Plugging the expansions into \eqref{eqn:ch1} and \eqref{eqn:chem} and matching the coefficients of $\e$ powers in both equations, the $O(\frac{1}{\e^2})$ of \eqref{eqn:ch1}  and \eqref{eqn:chem} yields
\begin{eqnarray}
\label{eqn:u0}&&0=\nabla \cdot  \left(\nabla \mu^{(0)}+w^{(0)}\nabla u^{(0)}\right)\\
\label{eqn:mu0}&&\mu^{(0)}=q'(u^{(0)}).
\end{eqnarray}
Since
$$
 w^{(0)}=\mu^{(0)}\frac{ g'(u^{(0)})}{g(u^{(0)})},
$$
then $u^{(0)}=1$ or $u^{(0)}=-1$ satisfies equations \eqref{eqn:u0}-\eqref{eqn:mu0}. In particular, such choice of $u^{(0)}$ implies $\mu^{(0)}=0$.

The $O(\frac{1}{\e})$ equation of \eqref{eqn:ch1}  and \eqref{eqn:chem} reduces to

\begin{eqnarray}
\label{eqn:u1}&&0=\nabla \cdot  \left(\nabla \mu^{(1)}+w^{(0)}\nabla u^{(1)}+w^{(1)}\nabla u^{(0)}\right),\\
\label{eqn:mu1}&&\mu^{(1)}=q''(u^{(0)})u^{(1)}.
\end{eqnarray}
Since $u^{(0)}=1$ or $-1$,
 $u^{(1)}=0$ satisfies \eqref{eqn:u1}-\eqref{eqn:mu1}. Moreover, such choice of $u^{(1)}$ guarantees $\mu^{(1)}=0$.

The $O(1)$ equation of \eqref{eqn:ch1}  and \eqref{eqn:chem}, taking into account of the fact $u^{(0)}=\pm 1$, $\mu^{(0)}=\mu^{(1)}=0$, reduces to

\begin{eqnarray}
\label{eqn:u2}&&0=\nabla \cdot  \left(\nabla \mu^{(2)}+w^{(0)}\nabla u^{(2)}\right),\\
\label{eqn:mu2}&&\mu^{(2)}=q''(u^{(0)})u^{(2)}.
\end{eqnarray}

Thus $u^{(2)}=0$ satisfies  \eqref{eqn:u2}-\eqref{eqn:mu2}. Moreover, such choice of $u^{(1)}$ guarantees $\mu^{(2)}=0$.

In general, if $u^{(0)}=\pm 1$, $u^{(1)}=u^{(2)}=\hdots=u^{(k+1)}=0$, the $O(\e^k)$ for $k\geq 1$ equation of \eqref{eqn:ch1}  and \eqref{eqn:chem} yields
\begin{eqnarray}
\label{eqn:uk} && 0=\nabla \cdot  \left(\nabla \mu^{(k+2)}+w^{(0)}\nabla u^{(k+2)}\right)\\
\label{eqn:muk} && \mu^{(k+2)}=q''(u^{(0)})u^{(k+2)}.
\end{eqnarray}
Thus $u^{(k+2)}=0$ satisfies \eqref{eqn:uk} and\eqref{eqn:muk}.

In summary, the $u=1$ or $u=-1$ in the outer region.

\subsection{Inner expansions}
For the small inner regions near the dislocations, we introduce local coordinates near the dislocations. Considering a dislocation $C$ parameterized by arc length parameter $s$. We denote a point on the dislocation by $\rr_0(s)$ with  tangent unit vector $\mathbf{t}(s)$ and inward normal vector $\bn(s)$. A point near the dislocation $C$ is expressed as
$$
\rr(s,d)=\rr_0(s)+d\mathbf{n}(s),
$$\
where $d$ is the signed distance from point $\rr$ to the dislocation. Since the gradients fields are of order $O(\frac{1}{\e})$, we introduce $\rho=\frac{d}{\e}$ and use coordinates $(s,\rho)$ in the inner region.
Under this setting, we write $u(x,y,t)=U(s,\rho,t)$ and equation \eqref{eqn:ch} can be written as
\begin{eqnarray}
&&\label {eqn:chin}\hspace{0.5 in} g(U)\left(\partial_t U+\frac{1}{\e}v_n\partial_{\rho}U\right)=\frac{M_0}{1-\e\rho\kappa}\partial_s\left(\frac{1}{1-\e\rho\kappa}\left(\partial_s \mu-\mu\frac{g'(U)}{g(U)}\partial_s U\right)\right) \\
&& \notag \hspace{1.9 in}+\frac{1}{\e^2}\frac{M_0}{1-\e\rho\kappa}\partial_{\rho}\left((1-\e\rho\kappa)\left(\partial_{\rho}\mu-\mu \frac{g'(U)}{g(U)}\partial_{\rho}U\right)\right),\\
&&\label{eqn:chemin} \hspace{1 in} \mu=-\frac{1}{1-\e\rho\kappa}\partial_s\left(\frac{1}{1-\e\rho\kappa}\partial_sU\right)-\frac{1}{\e^2}\frac{1}{1-\e\rho\kappa}\partial_{\rho}\left((1-\e\rho\kappa)\partial_{\rho}U\right)+\frac{1}{\e^2}q'(U).
\end{eqnarray}

Assume $\mu$ takes the same form expansion as  \eqref{eqn:muexp} and the following expansions hold for $U$  within dislocation core region:
\begin{equation}
U(s,\rho,t)=U^{(0)}(\rho)+\e U^{(1)}(s,\rho,t)+\e^2 U^{(2)}+\hdots. \label{eqn:U}
\end{equation}

Here we assume the leading order solution $U^{(0)}$, which describes the dislocation core profile, remains the same at all points on the dislocation at any time.

Set $$W=\mu\frac{g'(U)}{g(U)}=\frac{1}{\e^2}\left(W^{(0)}+W^{(1)}\e+W^{(2)}\e^2+\hdots \right),$$ the leading order for equation \eqref{eqn:chin} and \eqref{eqn:chemin} is $O(\frac{1}{\e^4})$, which yields
\begin{eqnarray}
&&\label{eqn:U0} 0= \partial_{\rho}\left(\partial_{\rho}\mu^{(0)}-W^{(0)}\partial_{\rho}U^{(0)}\right),\\
&& \label{eqn:mu0in} \mu^{(0)}=-\partial_{\rho
\rho}U^{(0)}+q'(U^{(0)}).
\end{eqnarray}
Substituting $W^{(0)}=\mu^{(0)}\frac{g'(U^{(0)})}{g(U^{(0)})}$ into \eqref{eqn:U0}, we can rewrite \eqref{eqn:U0} as $$0=\partial_{\rho}\left(\partial_{\rho}\mu^{(0)}-\mu^{(0)}\partial_{\rho}\ln g(U^{(0)})\right).$$
Integrating this equation, we have
\begin{equation}
\partial_{\rho}\mu^{(0)}-\mu^{(0)}\partial_{\rho}\ln g(U^{(0)})=C_0(s). \label{eqn:c0}
\end{equation}
Since $\mu^{(0)}=0$ in the outer region, we must have $\mu^{(0)}\rightarrow 0$ and $\partial_{\rho}\mu^{(0)}\rightarrow 0$ as $\rho\rightarrow \pm \infty$. Therefore $C_0(s)=0$. Dividing \eqref{eqn:c0} by $\mu^{(0)}$ and integrating, we have $\mu^{(0)}=\tilde C_0(s)g(U^{(0)})$. Since $\mu^{(0)}$ is independent of $s$ and is $0$ in the outer region, we must have $\tilde C_0(s)=0$. Thus
\begin{equation}
\mu^{(0)}=-\partial_{\rho\rho}U^{(0)}+q'(U^{(0)})=0. \label{eqn:Umu0}
\end{equation}
Solution $U^{(0)}$  to \eqref{eqn:Umu0} subject to far field condition $U^{(0)}{(+\infty)}=-1$ and $U^{(0)}{(-\infty)}=1$ can be found numerically (see \cite{NXY20} for example). In particular, $\partial_{\rho}U^{(0)}<0$ for all $\rho$.

Next, the $O(\frac{1}{\e^3})$ equation of \eqref{eqn:chin} and \eqref{eqn:chemin} yields, using $\mu^{(0)}=0$, that
\begin{eqnarray}
&&\label{eqn:U1} 0= \partial_{\rho}\left(\partial_{\rho}\mu^{(1)}-W^{(1)}\partial_{\rho}U^{(0)}\right),\\
&& \label{eqn:mu1in} \mu^{(1)}=-\partial_{\rho\rho}U^{(1)}+\kappa \partial_{\rho}U^{(0)}+q''(U^{(0)})U^{(1)}.
\end{eqnarray}
When $\mu^{(0)}=0$, we have  $W^{(1)}=\mu^{(1)}\frac{g'(U^{(0)})}{g(U^{(0)})}$. Sustituting into \eqref{eqn:U1} and integrating, we have

\begin{equation}
\partial_{\rho}\mu^{(1)}-\mu^{(1)}\partial_{\rho}\ln g(U^{(0)})=C_1(s). \label{eqn:c1}
\end{equation}

Matching with the outer solutions ($ \partial_{\rho}\mu^{(1)}, \mu^{(1)} \rightarrow 0 $ as $\rho \rightarrow \pm \infty$), we conclude that $C_1(s)=0$. Dividing  \eqref{eqn:c1} by $\mu^{(1)}$ and integrating, we have $\mu^{(1)}=\tilde C_1(s)g(U^{(0)})$. Thus \eqref{eqn:mu1in}  can be written as
\begin{equation}
LU^{(1)}=-\kappa \partial_{\rho}U^{(0)}+\tilde C_1(s)g(U^{(0)})\label{eqn:L}
\end{equation}
where $L=-\partial_{\rho\rho}+q''(U^{(0)})$ is a linear operator whose  kernal is span$\{\partial_{\rho}U^{(0)}\}$.  \eqref{eqn:L} is sovlable iff the right hand side is perpendicular to the kernal of $L$, i.e.
$$
\int_{-\infty}^{+\infty}\left(-\kappa \partial_{\rho}U^{(0)}+\tilde C_1(s)g(U^{(0)})\right)\partial_{\rho}U^{(0)} d\rho=0.
$$
From this, we conclude
$$
\tilde C_1(s)=\alpha \kappa,
$$
where positive constants $\alpha$ and $\beta$ are given by
$$
\alpha=\frac{\int_{-\infty}^{+\infty}\left(\partial_{\rho}U^{(0)}\right)^2 d\rho}{\int_{-\infty}^{+\infty}g(U^{(0)})\partial_{\rho}U^{(0)}d\rho}<0.
$$
Therefore
\begin{equation}
\mu^{(1)}=\alpha \kappa g(U^{(0)}). \label{eqn:mu1exp}
\end{equation}
Letting $\overline{\mu}=\frac{\mu}{g(U)}$, \eqref{eqn:chin} can be written as
\begin{eqnarray}\notag
&& g(U)\left(\partial_t U+\frac{1}{\e}v_n\partial_{\rho}U\right)\\ \label{eqn:chinmubar}
&=& \frac{M_0}{1-\e\rho\kappa}\partial_s\left(\frac{g(U)}{1-\e\rho\kappa}\left(\partial_s \overline \mu\right) \right)+\frac{1}{\e^2}\frac{M_0}{1-\e\rho\kappa}\partial_{\rho}\left(\left(1-\e\rho\kappa\right)g(U)\partial_{\rho}\overline\mu\right)
\end{eqnarray}
Using $\mu^{(0)}=0$, $\partial_{\rho}\overline \mu^{(1)}=\partial_{\rho}\frac{\mu^{(1)}}{g(U^{(0)})}=0$, the $O(\frac{1}{\e^2})$ order equation of \eqref{eqn:chinmubar} reduces to
$$
\partial_{\rho}\left(g(U^{(0)})\partial_{\rho}U^{(2)}\right)=0.
$$
Integrating with respect to $\rho$, we have $g(U^{(0)})\partial_{\rho}U^{(2)}=C_2(s)$. Matching with outer solutions, we must have $C_2(s)=0$. Thus $\partial_{\rho}U^{(2)}=0$ which gives  $U^{(2)}=\tilde C_2(s)$.

Next we look at the $O(\frac{1}{\e})$ equation of \eqref{eqn:chinmubar}. Using $\mu^{(0)}=0$, $\partial_{\rho}\overline \mu^{(1)}=0$ and $\partial_{\rho}\overline \mu^{(2)}=0$, we have
$$
g(U^{(0)})v_n\partial_{\rho} U^{(0)}=M_0\partial_s\left(g(U^{(0)})\partial_s\overline \mu^{(1)}\right)+M_0\partial_{\rho}\left(g(U^{(0)})\partial_{\rho} \overline\mu^{(3)}\right).
$$
Integrating this equation with respect to $\rho$ and matching with outer solutions yields
\begin{equation}
v_n=\lambda \partial_{ss}\overline \mu^{(1)} \label{eqn:vn}
\end{equation}
where we used the fact that $g(U^{(0)})$ is independent of $s$ and
$$
\lambda=\frac{M_0\int_{-\infty}^{+\infty}g(U^{(0)})d\rho}{\int_{-\infty}^{+\infty}g(U^{(0)})\partial_{\rho}U^{(0)}d\rho}<0.
$$
By \eqref{eqn:mu1exp}, we have $\overline\mu^{(1)}=\alpha \kappa$, substitute this into \eqref{eqn:vn}, we obtain the sharp interface limit equation
\begin{equation}
v_n=\lambda \partial_{ss}\left(\alpha \kappa\right). \label{eqn:sharpif}
\end{equation}
\begin{remark}
Notice here the outer and inner expansions are similar to the expansions in \cite{NXY20}. We wrote out all details here for readers' convenience.
\end{remark}

\section {Weak solution for phase field model with positive mobilities}
\label{sec:ndg}
In this section we prove existence of weak solutions for phase field model with positive mobilities. Let $\mathbb Z_{+}$ be the set of nonnegative integers and $\Omega=[0,2\pi]^n$ with $n\leq 2$.  We pick an orthonormal basis for $L^2(\Omega)$ as
\begin{eqnarray*}
\{\phi_j:j=1,2,\hdots\}=\left \{(2\pi)^{-{n/2}},\text{Re}\left(\pi^{-n/2}e^{i\xi\cdot x}\right), \text{Im}\left(\pi^{-n/2}e^{i\xi\cdot x}\right):\xi \in \mathbb Z^n_{+}\backslash \{0,\hdots,0\}\right \}.
\end{eqnarray*}
Observe $\{\phi_j\}$ is also orthogonal in $H^k(\Omega)$ for any $k\geq 1$.  Here and throughout the paper, we denote  $\Omega_T= (0,T)\times \Omega$.
\subsection{Galerkin approximations}

Define
$$
u^N(x,t)=\sum_{j=1}^Nc_j^N(t)\phi_j(x), \hspace{0.5 in} \mu^N(x,t)=\sum_{j=1}^Nd^N_{j}(t)\phi_j(x),
$$
where $\{c_j^N,d^N_j\}$ satisfy
\begin{eqnarray}\label{eqn:uN}
\int_{\Omega}\partial_t u^N\phi_j dx&=&-\int_{\Omega}M_{\theta}(u^N)\nabla \frac{\mu^N}{g_{\theta}(u^N)}\cdot \nabla \frac{\phi_j}{g_{\theta}(u^N)}dx,\\ \label{eqn:muN}
\int_{\Omega}\mu^N\phi_j dx&=&\int_{\Omega}\left(\nabla u^N\cdot \nabla \phi_j+q'(u^N)\phi_j\right)dx,\\ \label{eqn:uN0}
u^N(x,0)&=&\sum_{j=1}^N\left(\int_{\Omega} u_0\phi_j dx\right)\phi_j(x).
\end{eqnarray}
\eqref{eqn:uN}-\eqref{eqn:uN0} is an initial value problem for a system of ordinary equations for $\{c_j^N(t)\}$. Since right hand side of \eqref{eqn:uN} is continuous in $c_j^N$, the system has a local solution.

Define energy functional
$$
E(u)=\int_{\Omega}\left\{\frac{1}{2}|\nabla u|^2 +q(u)\right\} dx.
$$
Direct calculation yields
$$
\frac{d}{dt}E(u^N(x,t))=-\int_{\Omega}M_{\theta}(u^N)\left|\nabla \frac{\mu^N}{g_{\theta}(u^N)}\right|^2 dx,
$$
integration over $t$ gives the following energy identity.
\begin{eqnarray}\notag
&&\int_{\Omega}\left(\frac{1}{2}|\nabla u^N(x,t)|^2+q(u^N(x,t))\right)dx\\ \notag
&&+\int_0^t\int_{\Omega}M_{\theta}(u^N(x,\tau))\left|\nabla \frac{\mu^N(x,\tau)}{g_{\theta}(u^N(x,\tau))}\right|^2 dxd\tau\\ \label{eqn:energyid}
&=&\int_{\Omega}\left(\frac{1}{2}|\nabla u^N(x,0)|^2+q(u^N(x,0))\right)dx\\ \notag
&\leq&\left \Arrowvert \nabla u_0\right \Arrowvert_{L^2(\Omega)}^2+C\left(\left\Arrowvert u_0\right\Arrowvert^{r+1}_{H^1{\Omega}}+|\Omega|\right)\leq C <\infty
\end{eqnarray}
Here and throughout the paper, $C$ represents a generic constant possibly depending only on  $T$, $\Omega$,  $u_0$ but not on $\theta$. Since $\Omega$ is bounded region, by  growth assumption assumption  \eqref{eqn:qu} and Poincare's inequality, the energy identity \eqref{eqn:energyid} implies
$
u^N\in L^{\infty}(0,T;H^1(\Omega))
$
with
\begin{equation}
\left\Arrowvert u^N\right\Arrowvert _{L^{\infty}(0,T;H^1(\Omega))}\leq C \text{ for all } N, \label{eqn:uNbd}
\end{equation}
and
\begin{equation}
\left \Arrowvert \sqrt{M_{\theta}(u^N)}\nabla \frac{\mu^N}{g_{\theta}(u^N)}\right \Arrowvert_{L^2(\Omega_T)}\leq C \text{ for all } N.  \label{eqn:muNbd}
\end{equation}
By \eqref{eqn:uNbd}, the coefficients $\{c_j^N(t)\}$ are bounded in time, thus the system \eqref{eqn:uN}-\eqref{eqn:uN0} has a global solution. In addition, by Sobolev embedding theorem and growth assumption \eqref{eqn:q'u} on $q'(u)$, we have
$$
q'(u^N)\in L^{\infty}(0,T; L^p(\Omega)), \hspace{0.5in} M_{\theta}(u^N)\in L^{\infty}(0,T;L^p(\Omega))
$$
 for any $1\leq p<\infty$ with
\begin{eqnarray} \label{eqn:qunbd}
&&\left \Arrowvert q'(u^N)\right\Arrowvert_{L^{\infty}(0,T;L^p(\Omega))}\leq C  \text{ for all } N, \\ \label{eqn:mbd}
&&\left \Arrowvert M_{\theta}(u^N)\right \Arrowvert_{L^{\infty}(0,T;L^p(\Omega))}\leq C     \text{ for all } N.
\end{eqnarray}
\subsection{Convergence of $u^N$}\label{sec:uNconv}

Given $q>2$ and any $\phi\in L^2(0,T;W^{1,q}(\Omega))$, let  $\Pi_N\phi(x,t)=\sum_{j=1}^N\left(\int_{\Omega}\phi(x,t)\phi_j(x)dx\right) \phi_j(x)$ be the orthogonal projection of $\phi$ onto span$\{\phi_j\}_{j=1}^N$.  Then
\begin{eqnarray*}
&&\left \arrowvert \int_{\Omega}\partial_tu^N\phi dx\right \arrowvert=\left\arrowvert \int_{\Omega}\partial_tu^N\Pi_N\phi dx\right \arrowvert=\left \arrowvert \int_{\Omega}M_{\theta}(u^N)\nabla \frac{\mu^N}{g_{\theta}(u^N)}\cdot \nabla \frac{\Pi_N\phi}{g_{\theta}(u^N)}dx\right \arrowvert\\
&\leq& \left(\int_{\Omega}M_{\theta}(u^N)\left\arrowvert \nabla \frac{\mu^N}{g_{\theta}(u^N)}\right\arrowvert^2dx\right)^{\frac{1}{2}}\left(\int_{\Omega}M_{\theta}(u^N)\left \arrowvert \nabla\frac{\Pi_N\phi}{g_{\theta}(u^N)}\right \arrowvert^2 dx\right)^{\frac{1}{2}}.
\end{eqnarray*}

Since
$$
 \nabla\frac{\Pi_N\phi}{g_{\theta}(u^N)}=\frac{1}{g_{\theta}(u^N)}\nabla \Pi_N \phi-\Pi_N \phi \frac{g_{\theta}'(u^N)}{g^2_{\theta}(u^N)}\nabla u^N,
$$
we have
\begin{eqnarray*}
&& \int_{\Omega}M_{\theta}(u^N)\left \arrowvert \nabla\frac{\Pi_N\phi}{g_{\theta}(u^N)}\right \arrowvert^2 dx\\
&\leq& 2 M_0\int_{\Omega}\left(\frac{1}{g_{\theta}(u^N)}\left \arrowvert \nabla \Pi_N \phi\right \arrowvert ^2 + \frac{|g'_{\theta}(u^N)|^2}{g^3_{\theta}(u^N)} | \Pi_N \phi|^2|\nabla u^N|^2\right)dx\\
&\leq& C(M_0,\theta)\left(\left \Arrowvert \nabla \Pi_N \phi \right \Arrowvert^2_{L^2(\Omega)}+ \left \Arrowvert \Pi_N \phi \right\Arrowvert^2_{L^{\infty}(\Omega)}\left \Arrowvert \nabla u^N \right \Arrowvert^2_{L^2(\Omega)}\right)  \text{ here is where we need }  q>2\\
&\leq& C(M_0,\theta)\left(\left \Arrowvert \Pi_N \phi \right \Arrowvert^2_{W^{1,q}(\Omega)}\right) \leq C(M_0,\theta)\left \Arrowvert\phi \right \Arrowvert^2_{W^{1,q}(\Omega)}.
\end{eqnarray*}
Therefore
\begin{equation}
\left \Arrowvert \partial_t u^N \right \Arrowvert_{L^2(0,T; (W^{1,q}(\Omega))')} \leq C(M_0,\theta) \text{ for all } N. \label{eqn:uNtbd}
\end{equation}

 For $1\leq s < \infty$, since $n\leq 2$, by Sobolev embedding theorem and Aubin-Lions Lemma (see \cite{Sim86} and  Remark \ref{rmk:cpt}) , the following embeddings are compact :
$$
\left \{f\in L^2(0,T;H^1(\Omega)):\partial_tf\in L^2(0,T;(W^{1,q}(\Omega))'\right \}\hookrightarrow  L^2(0,T;L^s(\Omega)),
$$
and
$$
\left \{f\in L^{\infty}(0,T;H^1(\Omega)):\partial_tf\in L^2(0,T;(W^{1,q}(\Omega))'\right \}\hookrightarrow  C([0,T];L^s(\Omega)).
$$
From this  and  the boundedness of $\{u^N\}$ and $\{\partial_t u^N\}$, we can find a subsequence, and $u_{\theta} \in L^{\infty}(0,T;H^1(\Omega))$ such that as $N\rightarrow \infty$,
\begin{eqnarray}\label{eqn:uNwC}
u^N &\rightharpoonup& u_{\theta}\text{  weakly-* in }L^{\infty}(0,T;H^1(\Omega)), \\ \label{eqn:uNC}
u^N &\rightarrow& u_{\theta} \text{ strongly in } C([0,T];L^s(\Omega)),\\ \label{eqn:uNp}
u^N&\rightarrow &
u_{\theta} \text{ strongly in } L^2(0,T;L^s(\Omega)) \text{ and  a.e.  in }\Omega_T, \\\label{eqn:uNtC}
\partial_tu^N &\rightharpoonup &\partial_t u_{\theta} \text{ weakly in } L^2(0,T;(W^{1,q}(\Omega))')
\end{eqnarray}
for $1\leq s<\infty$. In addition
$$
\left \Arrowvert u_{\theta} \right \Arrowvert_{L^{\infty}(0,T;H^1(\Omega))} \leq C, \hspace{0.5 in} \left \Arrowvert \partial_t u_{\theta}\right \Arrowvert_{L^2(0,T; (W^{1,q}(\Omega))')} \leq C(M_0, \theta).
$$

By \eqref{eqn:uNC}, growth assumption \eqref{eqn:q'u} on $q'(u^N)$, and general dominated convergence Theorem, we have for any  $1\leq s<\infty$,
\begin{eqnarray}\label{eqn:MuNC}
&& M_{\theta}(u^N) \rightarrow M_{\theta}(u_{\theta}) \text{ strongly in } C([0,T]; L^s(\Omega)),\\ \label{eqn:sqMC}
&& \sqrt {M_{\theta}(u^N)} \rightarrow \sqrt {M_{\theta}(u_{\theta})} \text{ strongly in } C([0,T]; L^s(\Omega)),\\ \label{eqn:quNC}
&& q'(u^N) \rightarrow q'(u_{\theta}) \text{ strongly in } C([0,T];L^s(\Omega)).
\end{eqnarray}
 By \eqref{eqn:qunbd} and \eqref{eqn:quNC}, we have
\begin{equation}
q'(u^N) \rightharpoonup q'(u_{\theta}) \text{ weakly-* in } L^{\infty}([0,T];L^s(\Omega)). \label{eqn:quNwC}
\end{equation}
\begin{remark}\label{rmk:cpt}
Let $X$, $Y$, $Z$ be Banach spaces with compact embedding $X\hookrightarrow Y$ and continuous embedding $Y\hookrightarrow Z$. Then the embeddings
\begin{equation}
\{f \in L^p(0,T;X); \partial_t f \in L^1(0,T;Z)\} \hookrightarrow L^p(0,T;Y) \label{eqn:l2embed}
\end{equation}
and
\begin{equation}
\{f \in L^{\infty}(0,T;X); \partial_t f \in L^r(0,T;Z)\} \hookrightarrow C([0,T];Y) \label{eqn:linftyembed}
\end{equation}
are compact for any $1\leq p < \infty$ and $r>1$ (Corollary 4, \cite{ Sim86}, see also \cite{Lions69}) . For convergence of $u^N$, we apply this for $p=2=r$ with  $X=H^1(\Omega)$, $Y=L^s(\Omega)$ for $1\leq s<\infty$ and $Z=W^{1,q}(\Omega)'$.
\end{remark}
\subsection{Weak solution }

By  \eqref{eqn:muNbd} and the lower bound on $M_{\theta}$, we have
\begin{equation}
\left \Arrowvert \nabla \frac{\mu^N}{g_{\theta}(u^N)}\right \Arrowvert_{L^2(\Omega_T)} \leq C \theta^{-\frac{m}{2}}. \label{eqn:gramuNg}
\end{equation}

By \eqref{eqn:muN}, \eqref{eqn:uNbd} and \eqref{eqn:qunbd}, we have
\begin{eqnarray}\notag
&& \left \arrowvert \int_{\Omega}\frac{\mu^N\phi_1}{g_{\theta}(u_N)}\right \arrowvert dx=\left \arrowvert \int_{\Omega}\mu^N \Pi_N\left(\frac{\phi_1}{g_{\theta}(u^N)}\right)\right \arrowvert dx\\ \label{eqn:muNg}
&=&\left \arrowvert \int_{\Omega}\nabla u^N\cdot \nabla \Pi_N\left (\frac{\phi_1}{g_{\theta}(u^N)}\right)dx+\int_{\Omega}q'(u^N)\Pi_N\left(\frac{\phi_1}{g_{\theta}(u^N)}\right)dx\right \arrowvert \\ \notag
&=&\left\arrowvert\int_{\Omega}\nabla u^N\cdot \nabla \left (\frac{\phi_1}{g_{\theta}(u^N)}\right)dx+\int_{\Omega}q'(u^N)\Pi_N\left(\frac{\phi_1}{g_{\theta}(u^N)}\right)dx\right \arrowvert \\ \notag
&\leq& C\left(\theta^{-m-1}\left \Arrowvert \nabla u^N\right \Arrowvert_{L^2(\Omega)}^2+\theta^{-m}\left \Arrowvert q'(u^N)\right\Arrowvert_{L^2(\Omega)}\left \Arrowvert \phi_1\right \Arrowvert_{L^2(\Omega)}\right)\\ \notag
&\leq& C\theta^{-m-1}.
\end{eqnarray}
 \eqref{eqn:gramuNg},\eqref{eqn:muNg} and Poincare's inequality yield
$$
\left \Arrowvert \frac{\mu^N}{g_{\theta}(u^N)} \right \Arrowvert_{L^2(0,T;H^1(\Omega))} \leq C(\theta^{-m-1}+1).
$$
Thus there exists a $w_{\theta} \in L^2(0,T;H^1(\Omega))$ and a subsequence of $\frac{\mu^N}{g_{\theta}(u^N)}$, not relabeled, such that
\begin{equation}
\frac{\mu^N}{g_{\theta}(u^N)}\rightharpoonup w_{\theta} \text{ weakly in } L^2(0,T;H^1(\Omega)). \label{eqn:mugNC}
\end{equation}
Therefore by \eqref{eqn:MuNC}, \eqref{eqn:mugNC} and Sobolev embedding theorem, we have
\begin{equation}
\mu^N=g_{\theta}(u^N)\cdot \frac{\mu^N}{g_{\theta}(u^N)}\rightharpoonup \mu_{\theta}=g_{\theta}(u_\theta)w_{\theta} \text{ weakly in } L^2(0,T;W^{1,s}(\Omega))\label{eqn:muNC}
\end{equation}
for any $1\leq s<2$. Combining \eqref{eqn:sqMC}, \eqref{eqn:mugNC}and  \eqref{eqn:muNC}, we have
\begin{equation}
\sqrt{M_{\theta}(u^N)}\nabla \frac{\mu^N}{g_{\theta}(u^N)} \rightharpoonup \sqrt{M_{\theta}(u_\theta)}\nabla \frac{\mu_\theta}{g_{\theta}(u_\theta)} \text{ weakly in } L^2(0,T;L^q(\Omega)) \label{eqn:nonlinearwC}
\end{equation}
for any $1\leq q <2$. By \eqref{eqn:muNbd}, we can improve this convergence to
\begin{equation}
\sqrt{M_{\theta}(u^N)}\nabla \frac{\mu^N}{g_{\theta}(u^N)} \rightharpoonup \sqrt{M_{\theta}(u_\theta)}\nabla \frac{\mu_\theta}{g_{\theta}(u_\theta)} \text{ weakly in } L^2(0,T;L^2(\Omega)). \label{eqn:muNl2wC}
\end{equation}
By \eqref{eqn:muN}, we have
$$
\int_{\Omega} \mu^N u^N dx=\int_{\Omega}\left(|\nabla u^N|^2 dx +q'(u^N)u^N\right)dx.
$$
Integrating with respect to $t$ from $0$ to $T$, we have on $\Omega_T=\Omega\times [0,T]$,
\begin{eqnarray*}
\int_{\Omega_T} \mu^N(x,\tau) u^N(x,\tau) dxd\tau
=\int_{\Omega_T}\left(\nabla u^N(x,\tau)|^2 dx +q'(u^N(x,\tau))u^N(x,\tau)\right)dxd\tau.
\end{eqnarray*}
Passing to the limit in the equation above, by \eqref{eqn:uNp}, \eqref{eqn:quNC} and \eqref{eqn:muNC}, we have
\begin{eqnarray}\label{eqn:guNC}
\int_{\Omega_T} \mu_{\theta}u_{\theta} dxd\tau=\lim_{N\rightarrow \infty}\int_{\Omega_T}|\nabla u^N|^2 dxd\tau+\int_{\Omega_T}q'(u_{\theta})u_{\theta} dxd\tau
\end{eqnarray}
On the other hand,
\begin{eqnarray}\label{eqn:muNuN}
&&\hspace {0.2 in}\int_{\Omega_T}\mu^N(x,\tau) u_{\theta}(x,\tau)dxd\tau=\int_{\Omega_T}\mu^N(x,\tau) \Pi_Nu_{\theta}(x,\tau)dxd\tau \\ \notag
&=&\int_{\Omega_T}\left(\nabla u^N\cdot \nabla\Pi_N u_{\theta}(x,\tau)  +q'(u^N) \Pi_N u_{\theta}(x,\tau)\right) dxd \tau \\ \notag
&=&\int_{\Omega_T}\left(\nabla u^N\cdot \nabla u_{\theta}(x,\tau)  +q'(u^N) \Pi_N u_{\theta}(x,\tau)\right)dxd \tau.
\end{eqnarray}
Since $\Pi_Nu_{\theta} \rightarrow u_{\theta}$ strongly in $L^2(\Omega_T)$, by \eqref{eqn:uNwC},\eqref{eqn:quNwC} and \eqref{eqn:muNC}, as $N\rightarrow \infty$, \eqref{eqn:muNuN} yields
\begin{equation}
\int_{\Omega_T}\mu_{\theta}u_{\theta}dxd\tau=\int_{\Omega_T}\left(|\nabla u_{\theta}|^2+q'(u_{\theta}))u_{\theta}\right)dxd\tau. \label{eqn:eid}
\end{equation}
\eqref{eqn:guNC} and \eqref{eqn:eid} gives
\begin{equation}
\lim_{N\rightarrow \infty}\int_{\Omega_T}|\nabla u^N|^2 dxd\tau=\int_{\Omega_T}|\nabla u_{\theta}|^2 dxd\tau. \label{eqn:normC}
\end{equation}
 By \eqref{eqn:uNbd}, $\nabla u^N \rightharpoonup \nabla u_{\theta}$ weakly in $L^2(\Omega_T)$, thus \eqref{eqn:normC} implies
\begin{equation}
\nabla u^N \rightarrow \nabla u_{\theta}\text{ strongly in } L^2(\Omega_T).\label{eqn:guNsC}
\end{equation}

Since $g_{\theta}\geq \theta^{m}$,  \eqref{eqn:uNp} implies
\begin{equation}
\frac{g'(u^N)}{g_{\theta}^{\frac{3}{2}}(u^N)}\rightarrow \frac{g'_{\theta}(u_{\theta})}{g_{\theta}^{\frac{3}{2}}(u_{\theta})} \text{ a.e  in  } \Omega_T. \label{eqn:lngC}
\end{equation}
In addition,  $\frac{g'(u^N)}{g_{\theta}^{\frac{3}{2}}(u^N)}$ is bounded by
\begin{equation}
\left \arrowvert \frac{g'(u^N)}{g_{\theta}^{\frac{3}{2}}(u^N)}\right \arrowvert \leq C \theta^{-1-\frac{m}{2}}. \label{eqn:lngbd}
\end{equation}

It follows from \eqref{eqn:guNsC}, \eqref{eqn:lngC}, \eqref{eqn:lngbd}  and generalized dominated convergence theorem (see Remark \ref{rmk:gdc}) that
\begin{equation}
\frac{g'(u^N)}{g_{\theta}^{\frac{3}{2}}(u^N)}\nabla u^N\rightarrow \frac{g'_{\theta}(u_{\theta})}{g_{\theta}^{\frac{3}{2}}(u_{\theta})}\nabla u_{\theta} \text{ strongly in } L^2(\Omega_T).  \label{eqn:l2conv}
\end{equation}
Let
$$
f^N(t)=\left\Arrowvert\frac{g'(u^N(x,t))}{g_{\theta}^{\frac{3}{2}}(u^N(x,t))}\nabla u^N(x,t)- \frac{g'_{\theta}(u_{\theta}(x,t))}{g_{\theta}^{\frac{3}{2}}(u_{\theta}(x,t))}\nabla u_{\theta}(x,t) \right \Arrowvert_{L^2(\Omega)},
$$
by \eqref{eqn:l2conv}, we can extract  a subsequence of $f^N$, not relabeled, such that  $f^N(t)\rightarrow 0$ a.e. in (0,T). By Egorov's theorem, for any given  $\delta>0$, there exists $T_{\delta}\subset [0,T]$ with $|T_{\delta}|<\delta$ such that $f^N(t)$ converges to $0$ uniformly on $[0,T]\backslash T_{\delta}$.

Given  $\alpha(t) \in L^2(0,T)$, for any $\varepsilon>0$, there exists $T_{\delta}\subset [0,T]$ with $|T_{\delta}|<\delta$ such that
\begin{equation}
\int_{T_\delta}\alpha^2(t) dt <\varepsilon. \label{eqn:alphat}
\end{equation}
Multiplying \eqref{eqn:uN} by $\alpha(t)$ and integrating  in time yield
\begin{eqnarray}  \label{eqn:uNl}
&&\int_0^T\alpha(t)\int_{\Omega}\partial_tu^N\phi_j dxdt=\int_0^T\alpha(t)\int_{\Omega}M_{\theta}(u^N)\nabla \frac{\mu^N}{g_{\theta}(u^N)}\cdot \nabla \frac{\phi_j}{g_{\theta}(u^N)}dxdt \\ \notag
&=&\int_{\Omega_T}M_0 \alpha(t) \nabla \frac{\mu^N}{g_{\theta}(u^N)}\cdot \nabla \phi_j dxdt \\ \notag
&&-\int_{\Omega_T}\alpha(t) \sqrt{M_0}\phi_j \frac{g'_{\theta}(u^N)}{g^{\frac{3}{2}}_{\theta}(u^N)}\nabla u^N\cdot \sqrt{M_{\theta}(u^N)}\nabla \frac{\mu^N}{g_{\theta}(u^N)} dxdt\\ \notag
&&=I-II.
\end{eqnarray}
Since $\alpha(t) \nabla \phi_j \in L^2(0,T;L^2(\Omega))$, by \eqref{eqn:mugNC} and \eqref{eqn:muNC}, we have
\begin{equation}
I=\int_{\Omega_T}M_0 \alpha(t) \nabla \frac{\mu^N}{g_{\theta}(u^N)}\cdot \nabla \phi_j dxdt \rightarrow \int_{\Omega_T}M_0 \alpha(t) \nabla \frac{\mu_{\theta}}{g_{\theta}(u_{\theta})}\cdot \nabla \phi_j dxdt. \label{eqn:Ilimit}
\end{equation}
To prove convergence on $II$, observe

\begin{eqnarray} \label{eqn:IIlimit}
&&\int_{\Omega_T}\alpha(t)\phi_j \frac{g'_{\theta}(u^N)}{g^{\frac{3}{2}}_{\theta}(u^N)}\nabla u^N \sqrt{M_{\theta}(u^N)}\nabla \frac{\mu^N}{g_{\theta}(u^N)}dxdt \\ \notag
&&- \int_{\Omega_T}\alpha(t)\phi_j \frac{g'_{\theta}(u_\theta)}{g^{\frac{3}{2}}_{\theta}(u_\theta)}\nabla u_\theta \sqrt{M_{\theta}(u_{\theta})}\nabla \frac{\mu_{\theta}}{g_{\theta}(u_{\theta})} dxdt\\ \notag
&=&\int_{\Omega_T}\alpha(t)\phi_j \left(\frac{g'_{\theta}(u^N)}{g^{\frac{3}{2}}_{\theta}(u^N)}\nabla u^N-\frac{g'_{\theta}(u_\theta)}{g^{\frac{3}{2}}_{\theta}(u_\theta)}\nabla u_\theta\right)\cdot \sqrt{M_{\theta}(u^N)}\nabla \frac{\mu^N}{g_{\theta}(u^N)} dxdt\\\notag
&&+\int_{\Omega_T}\alpha(t)\phi_j \frac{g'_{\theta}(u_\theta)}{g^{\frac{3}{2}}_{\theta}(u_\theta)}\nabla u_\theta\cdot \left( \sqrt{M_{\theta}(u^N)}\nabla \frac{\mu^N}{g_{\theta}(u^N)}-\sqrt{M_{\theta}(u_\theta)}\nabla \frac{\mu_\theta}{g_{\theta}(u_\theta)} \right)dxdt\\ \notag
&=&II_1+II_2
\end{eqnarray}
From bound
\begin{eqnarray*}
&&\int_{\Omega_T}\left\arrowvert \alpha(t)\phi_j \frac{g'_{\theta}(u_\theta)}{g^{\frac{3}{2}}_{\theta}(u_\theta)}\nabla u_\theta\right\arrowvert^2 dxdt \\
&\leq& C\theta^{-2-m} \left\Arrowvert \nabla u_{\theta}\right \Arrowvert_{L^{\infty}(0,T;L^2(\Omega))}^2 \int_0^T\alpha^2(t)^2 dt,
\end{eqnarray*}
 we conclude that $\alpha(t)\phi_j \frac{g'_{\theta}(u_\theta)}{g^{\frac{3}{2}}_{\theta}(u_\theta)}\nabla u_\theta \in L^2(\Omega_T)$. By \eqref{eqn:muNl2wC}, we can pass to the limit in $II_2$ and conclude
$$
II_2=\int_{\Omega_T}\alpha(t)\phi_j \frac{g'_{\theta}(u_\theta)}{g^{\frac{3}{2}}_{\theta}(u_\theta)}\nabla u_\theta\cdot \left( \sqrt{M_{\theta}(u^N)}\nabla \frac{\mu^N}{g_{\theta}(u^N)}-\sqrt{M_{\theta}(u_\theta)}\nabla \frac{\mu_\theta}{g_{\theta}(u_\theta)} \right)dxdt\rightarrow 0.
$$

To pass to the limit in $II_1$, we write
\begin{eqnarray*}
II_1&=&\int_{\Omega_T}\alpha(t)\phi_j \left(\frac{g'_{\theta}(u^N)}{g^{\frac{3}{2}}_{\theta}(u^N)}\nabla u^N-\frac{g'_{\theta}(u_\theta)}{g^{\frac{3}{2}}_{\theta}(u_\theta)}\nabla u_\theta\right)\cdot \sqrt{M_{\theta}(u^N)}\nabla \frac{\mu^N}{g_{\theta}(u^N)} dxdt\\
&=& \int_{T_{\delta}}\int_{\Omega}\alpha(t)\phi_j \left(\frac{g'_{\theta}(u^N)}{g^{\frac{3}{2}}_{\theta}(u^N)}\nabla u^N-\frac{g'_{\theta}(u_\theta)}{g^{\frac{3}{2}}_{\theta}(u_\theta)}\nabla u_\theta\right)\cdot \sqrt{M_{\theta}(u^N)}\nabla \frac{\mu^N}{g_{\theta}(u^N)} dxdt\\
&&+\int_{[0,T]\backslash T_{\delta}}\int_{\Omega}\alpha(t)\phi_j \left(\frac{g'_{\theta}(u^N)}{g^{\frac{3}{2}}_{\theta}(u^N)}\nabla u^N-\frac{g'_{\theta}(u_\theta)}{g^{\frac{3}{2}}_{\theta}(u_\theta)}\nabla u_\theta\right)\cdot \sqrt{M_{\theta}(u^N)}\nabla \frac{\mu^N}{g_{\theta}(u^N)} dxdt\\
&=&II_{11}+II_{12}.
\end{eqnarray*}

By \eqref{eqn:muNbd}, \eqref{eqn:uNwC}, \eqref{eqn:lngbd} and \eqref{eqn:alphat}, we can  bound $II_{11}$ by
\begin{eqnarray*}
|II_{11}|&\leq &\int_{T_{\delta}}|\alpha(t)| \left \Arrowvert\frac{g'_{\theta}(u^N)}{g^{\frac{3}{2}}_{\theta}(u^N)}\nabla u^N-\frac{g'_{\theta}(u_\theta)}{g^{\frac{3}{2}}_{\theta}(u_\theta)}\nabla u_\theta \right \Arrowvert_{L^2(\Omega)}\left\Arrowvert\sqrt{M_{\theta}(u^N)}\nabla \frac{\mu^N}{g_{\theta}(u^N)} \right\Arrowvert_{L^2(\Omega)} dt \\
&\leq& \left \Arrowvert \alpha(t)\right \Arrowvert_{L^2(T_{\delta})} \left \Arrowvert\sqrt{M_{\theta}(u^N)}\nabla \frac{\mu^N}{g_{\theta}(u^N)} \right \Arrowvert_{L^2(\Omega_T)}\left\Arrowvert \frac{g'_{\theta}(u^N)}{g^{\frac{3}{2}}_{\theta}(u^N)}\nabla u^N-\frac{g'_{\theta}(u_\theta)}{g^{\frac{3}{2}}_{\theta}(u_\theta)}\nabla u_\theta \right\Arrowvert_{L^{\infty}(0,T;L^2(\Omega))}\\
&\leq &C(\theta) \varepsilon.
\end{eqnarray*}
For $II_{12}$, we have
\begin{eqnarray*}
|II_{12}|&\leq &\int_{[0,T]\backslash T_{\delta}}|\alpha(t)| \left \Arrowvert\frac{g'_{\theta}(u^N)}{g^{\frac{3}{2}}_{\theta}(u^N)}\nabla u^N-\frac{g'_{\theta}(u_\theta)}{g^{\frac{3}{2}}_{\theta}(u_\theta)}\nabla u_\theta \right \Arrowvert_{L^2(\Omega)}\left\Arrowvert\sqrt{M_{\theta}(u^N)}\nabla \frac{\mu^N}{g_{\theta}(u^N)} \right\Arrowvert_{L^2(\Omega)} dt\\
&=&\int_{[0,T]\backslash T_{\delta}}|\alpha(t)|f^N(t)|\left\Arrowvert\sqrt{M_{\theta}(u^N)}\nabla \frac{\mu^N}{g_{\theta}(u^N)} \right\Arrowvert_{L^2(\Omega)} dt.
\end{eqnarray*}
Since  $f^N(t)$ converges to $0$ uniformly on $[0,T]\backslash T_{\delta}$, $\alpha(t) \in L^2(0,T)$ and $\left\Arrowvert\sqrt{M_{\theta}(u^N)}\nabla \frac{\mu^N}{g_{\theta}(u^N)} \right\Arrowvert_{L^2(\Omega_T)}\leq C $, letting $N \rightarrow \infty $ in $II_{12}$ yields $II_{12} \rightarrow 0$. Letting $\varepsilon \rightarrow 0$, we conclude $II_1\rightarrow 0$ as $N\rightarrow \infty$.
Passing to the limit in \eqref{eqn:uNl}, we have
\begin{eqnarray}\notag
&&\int_0^T\alpha (t) \int_{\Omega}\left <\partial_tu_{\theta},\phi_j \right>_{(W^{1,q}(\Omega))', W^{1,q}(\Omega))} dt\\ \label{eqn:phij}
&=&-\int_{\Omega_T}\alpha(t)M_{\theta}(u_{\theta})\nabla \frac{\mu_\theta}{g_{\theta}(u_\theta)}\cdot \nabla \frac{\phi_j}{g_\theta(u_\theta)} dxdt.
\end{eqnarray}
Fix $q>2$,  given any $\phi \in L^2(0,T;W^{1,q}(\Omega))$, its Fourier series $\sum_{j=1}^\infty a_j(t) \phi_j(x)$ converges strongly to $\phi$ in $L^2(0,T;W^{1,q}(\Omega))$. Hence
\begin{eqnarray}\label{eqn:erestimate}
&&\int_{\Omega_T}M_{\theta}(u_{\theta})\nabla \frac{\mu_{\theta}}{g_{\theta}(u_{\theta})}\cdot \nabla \frac{\phi-\Pi_N \phi}{g_{\theta}(u_{\theta})} dxdt \\ \notag
&=&\int_{\Omega_T}M_0\nabla \frac{\mu_{\theta}}{g_{\theta}(u_{\theta})}\cdot \nabla (\phi-\Pi_N \phi)dxdt \\ \notag
&&-\int_{\Omega_T}(\phi-\Pi_N \phi)\sqrt{M_0}\frac{g_{\theta}'(u_{\theta})}{g_{\theta}^{\frac{3}{2}}(u_{\theta})}\nabla u_{\theta}\cdot \sqrt{M_{\theta}(u_{\theta})}\nabla \frac{\mu_{\theta}}{g_{\theta}(u_{\theta})} dxdt\\ \notag
&=&J_1-J_2,
\end{eqnarray}
where $$
J_1=\int_{\Omega_T}M_0\nabla \frac{\mu_{\theta}}{g_{\theta}(u_{\theta})}\cdot \nabla (\phi-\Pi_N \phi)dxdt \rightarrow 0
$$ by\eqref{eqn:mugNC}, \eqref{eqn:muNC} and  strong convergence of $\Pi_N\phi $ to $\phi$ in $L^2(0,T;H^1(\Omega))$. We can bound $J_2$ by
\begin{eqnarray*}
&&|J_2|=\left \arrowvert \int_{\Omega_T}(\phi-\Pi_N \phi)\sqrt{M_0}\frac{g_{\theta}'(u_{\theta})}{g_{\theta}^{3/2}(u_{\theta})}\nabla u_{\theta}\cdot \sqrt{M_{\theta}(u_{\theta})}\nabla \frac{\mu_{\theta}}{g_{\theta}(u_{\theta})} dxdt\right \arrowvert  \\
&\leq &\sqrt{M_0}\int_0^T\left\Arrowvert\phi-\Pi_N \phi\right\Arrowvert_{L^{\infty}(\Omega)} \left\Arrowvert\frac{g_{\theta}'(u_{\theta})}{g_{\theta}^{3/2}(u_{\theta})}\nabla u_{\theta}\right\Arrowvert_{L^2(\Omega)} \left\Arrowvert\sqrt{M_{\theta}(u_{\theta})}\nabla \frac{\mu_{\theta}}{g_{\theta}(u_{\theta})}\right\Arrowvert_{L^2(\Omega)}\\
&\leq&\sqrt{M_0}\left\Arrowvert\frac{g_{\theta}'(u_{\theta})}{g_{\theta}^{3/2}(u_{\theta})}\nabla u_{\theta}\right\Arrowvert_{L^{\infty}(0,T;L^2(\Omega))}\left\Arrowvert\sqrt{M_{\theta}(u_{\theta})}\nabla \frac{\mu_{\theta}}{g_{\theta}(u_{\theta})}\right\Arrowvert_{L^2(\Omega_T)}\left\Arrowvert\phi-\Pi_N \phi\right\Arrowvert_{L^2(0,T;W^{1,q}(\Omega))}\\
&\rightarrow& 0 \text{ as } N\rightarrow \infty.
\end{eqnarray*}
Consequently \eqref{eqn:phij} and \eqref{eqn:erestimate} imply
\begin{equation}
\int_0^T\left <\partial_tu_{\theta},\phi \right>_{(W^{1,q}(\Omega))', W^{1,q}(\Omega))} dt=-\int_{\Omega_T}M_{\theta}(u_{\theta})\nabla \frac{\mu_\theta}{g_{\theta}(u_\theta)}\cdot \nabla \frac{\phi}{g_\theta(u_\theta)} dxdt \label{eqn:uthetaeq}
\end{equation}
for all $\phi \in L^2(0,T;W^{1,q}(\Omega))$ with $q>2$. Moreover, since $u^N(x,0)=\Pi_Nu_0(x) \rightarrow u_0(x)$ in $H^1(\Omega)$, we see that $u_\theta(x,0)=u_0(x)$ by \eqref{eqn:uNC}.

\begin{remark}\label{rmk:gdc}
(Generalized dominated convergence theorem) Assume  $E\subset \R^n$ is  measurable. $g_n\rightarrow g $ strongly in $L^q(E)$ for $1\leq q<\infty$ and $f_n$, $f$: $E\rightarrow \R^n$ are measurable functions satisfying
$$
f_n\rightarrow f \text{ a.e.  in } E; \hspace{0.1in} |f_n|^p \leq |g_n|^q \text{ a.e. in } E
$$
with $1\leq p< \infty$, then $f_n\rightarrow f$ in $L^p(E)$.
\end{remark}

\subsection{Regularity of $u_{\theta}$} \label{sec:regularity-u theta}
We now consider the regularity of $u_\theta$. Given any $a_j(t) \in L^2(0,T)$, $a_j(t)\phi_j \in L^2(0,T;C(\overline\Omega)$). Integrating   \eqref{eqn:muN} from $0$ to $T$, by \eqref{eqn:quNwC},\eqref{eqn:muNC} and  \eqref{eqn:guNsC}, we have
\begin{eqnarray*}
&&\int_{\Omega_T}\mu_\theta(x,t) a_j(t)\phi_j(x)dxdt
=\int_{\Omega_T}\left(\nabla u_\theta\cdot a_j(t)\nabla \phi_j+q'(u_\theta)a_j(t)\phi_j\right) dxdt
\end{eqnarray*}
for all $j \in \N$. Given any $\phi \in L^2(0,T;H^1(\Omega))$, its Fouirier series strongly converges to $\phi$ in $L^2(0,T;H^1(\Omega))$, therefore
\begin{eqnarray}\label{eqn:mutheta}
&&\int_{\Omega_T}\mu_\theta(x,t) \phi(x)dxdt
=\int_{\Omega_T}\left(\nabla u_\theta\cdot \nabla \phi+q'(u_\theta)\phi\right)dxdt.
\end{eqnarray}
Recall $\mu_{\theta} \in L^2(0,T;L^p(\Omega))$ and $q'(u_{\theta})\in L^{\infty}(0,T;L^p(\Omega))$ for any $1\leq p<\infty$, regularity theory implies  $u_\theta \in L^2(0,T;W^{2,p}(\Omega))$. Hence
\begin{equation}
\mu_\theta=-\Delta u_\theta+q'(u_\theta)\text{ a.e. in } \Omega_T.
\end{equation}
 Since growth assumption on $q$ implies $|q''(u)|\leq C(1+|u|^{r-1})$,  pick $p>2$, we have
\begin{eqnarray*}
&&\int_{\Omega}|\nabla q'(u_\theta)|^2 dx=\int_{\Omega}|q''(u_\theta)|^2|\nabla u_\theta|^2 dx \\
&\leq&\left \Arrowvert q''(u_\theta)\right \Arrowvert_{L^{\frac{2p}{p-2}}(\Omega)}^2\left \Arrowvert \nabla u_{\theta}\right \Arrowvert_{L^p(\Omega)}^2\\
&\leq &C\left(1+\left \Arrowvert u_\theta \right \Arrowvert_{L^{\frac{2p}{p-2}(r-1)}(\Omega)}^{2(r-1)}\right)\left \Arrowvert \nabla u_{\theta}\right \Arrowvert_{L^p(\Omega)}^2 \\
&\leq & C\left (1+\left \Arrowvert u_{\theta}\right \Arrowvert^{2(r-1)}_{L^{\infty}(0,T;H^1(\Omega))}\right )\left \Arrowvert \nabla u_{\theta}\right \Arrowvert_{L^p(\Omega)}^2
\end{eqnarray*}
Therefore $\nabla q'(u_\theta)=q''(u_\theta)\nabla u_\theta \in L^2(\Omega_T)$ with
\begin{eqnarray*}
\int_{\Omega_T}|\nabla q'(u_\theta)|^2 dxdt\leq \left (1+\left \Arrowvert u_{\theta}\right \Arrowvert^{2(r-1)}_{L^{\infty}(0,T;H^1(\Omega))}\right )\left \Arrowvert \nabla u_{\theta}\right \Arrowvert_{L^2(0,T;L^p(\Omega))}^2.
\end{eqnarray*}
Hence  $q'(u_\theta) \in L^2(0,T;H^1(\Omega))$, combined with $\mu_{\theta} \in L^2(0,T;W^{1,s}(\Omega))$ for any $1\leq s<2$, we have $u_\theta \in L^2(0,T;W^{3,s}(\Omega))$ and
\begin{equation}
\nabla \mu_\theta=-\nabla \Delta u_{\theta}+q''(u_\theta)\nabla u_\theta  \text{ a.e. in }  \Omega_T.
\end{equation}
Regularity of $u_{\theta}$ implies  $\nabla u_{\theta} \in L^{\infty}(0,T;L^2(\Omega))\cap L^2(0,T; L^{\infty}(\Omega))$. A simple interpolation shows  $\nabla u_{\theta} \in L^{\frac{2\mu}{\mu-2}}(0,T;L^{\mu}(\Omega))$ for any $\mu>2$.  Given any $\phi \in L^p(0,T;W^{1,q}(\Omega))$ with $p>2$ and $q>2$, we have  $g_{\theta}(u_{\theta})\phi \in L^2(0,T;W^{1.r}(\Omega))$ for any $r<\min(p,q)$. From this, we can pick  $g_{\theta}(u_{\theta})\phi$ as a test function in \eqref{eqn:uthetaeq}, we have
\begin{equation}
\int_{\Omega_T}\partial_tu_{\theta} g_{\theta}(u_\theta) \phi dxdt=-\int_{\Omega_T}M_{\theta}(u_\theta)\nabla \frac{\mu_\theta}{g_{\theta}(u_{\theta})}\cdot \nabla \phi dxdt \label{eqn:utheta2}
\end{equation}
for any $\phi \in L^p(0,T;W^{1,q}(\Omega))$ with $p,q>2$.
\begin{remark}
In fact, since $M_{\theta}(u_{\theta}) \in L^{\infty}(0,T;L^p(\Omega))$ for $1\leq p<\infty$, the right hand side of \eqref{eqn:utheta2} is well defined for any $\phi \in L^2(0,T,W^{1,q}(\Omega))$ for $q>2$ and we can extend  \eqref{eqn:utheta2} to hold for all $\phi \in L^2(0,T,W^{1,q}(\Omega))$.
\end{remark}
\subsection{Energy Inequality}

Since $u^N$ and $\mu^N$ satisfies energy identity \eqref{eqn:energyid}, passing to the limit as $N\rightarrow \infty$ and using strong convergence of $u^N(x,0)$ to $u_0$ in $H^1(\Omega)$, together with the weak convergence of $u^N$,
$q'(u^N)$ and $\sqrt{M_\theta(u^N) }\nabla \frac{\mu^N}{g_\theta(u^N)}$, the energy inequality \eqref{eqn:engineq} follows.

\section{Phase field model with degenerate mobility}
\label{sec:dg}

In this section, we prove theorem \ref{thm:dg}. Fix initial data $u_0 \in H^1(\Omega)$. We pick a montone decreasing  positive sequence $\theta_i$ with $\lim_{i\rightarrow \infty}\theta_i=0$. By theorem \ref{thm:ndg} and \eqref{eqn:utheta2}, for each $\theta_i$, there exists
$$
u_i \in L^{\infty}(0,T;H^1(\Omega))\cap L^2(0,T;W^{3,s}(\Omega))\cap C([0,T];L^p(\Omega))
$$
 with weak derivative
$$
\partial_t u_i \in L^2(0,T; (W^{1,q}(\Omega))'),
$$
where $1\leq p<\infty$, $1\leq s <2$, $q>2$ such that $u_{\theta_i}(x,0)=u_0(x)$ and for all $\phi \in L^2(0,T;W^{1,q}(\Omega))$,
\begin{eqnarray}
\label{eqn:ui}\int_0^T\int_{\Omega}\partial_t u_i\phi dxdt&=&-\int_0^T\int_{\Omega}M_i(u_i)\nabla \frac{\mu_i}{g_i(u_i)}\nabla \frac{\phi}{g_i(u_i)} dxdt,\\\mu_i&=&-\Delta u_i+q'(u_i).
\end{eqnarray}
Moreover,  for all $\psi \in L^p(0,T;W^{1,q}(\Omega))$ with $p,q >2$, the following holds:
\begin{equation} \label{eqn:ui2}
\int_0^T\int_{\Omega}g_i(u_i)\partial_t u_i\psi dxdt=-\int_0^T\int_{\Omega}M_i(u_i)\nabla \frac{\mu_i}{g_i(u_i)}\nabla \psi dxdt
\end{equation}
Here we write $u_i=u_{\theta_i}$, $M_i(u_i)=M_{\theta_i}(u_{\theta_i})$, $g_i(u_i)=g_{\theta_i}(u_{\theta_i})$ for simplicity of notations.

\subsection{Convergence of $u_i$ and equation for the limit function}
Noticing the bound in \eqref{eqn:uNbd} and \eqref{eqn:muNbd} only depends on $u_0$, we can find a constant $C$, independent of $\theta_i$ such that
\begin{eqnarray}\label{eqn:uibd}
&&\left \Arrowvert u_i\right \Arrowvert_{L^{\infty}(0,T;H^1(\Omega))}\leq C, \\ \label{eqn:Mui}
&& \left \Arrowvert \sqrt{M_i(u_i)}\nabla \frac{\mu_i}{g_i(u_i)}\right \Arrowvert_{L^2(\Omega_T)} \leq C.
\end{eqnarray}
Growth condition on $q'$,  and Sobolev embedding theorem give
\begin{eqnarray*}
&&\left \Arrowvert q'(u_i)\right \Arrowvert_{L^\infty(0,T;L^p(\Omega))} \leq C,\\
&&\left \Arrowvert M_i(u_i)\right \Arrowvert_{L^{\infty}(0,T;L^p(\Omega))} \leq C
\end{eqnarray*}
for any $1\leq p <\infty$. By \eqref{eqn:ui2}, for any $\phi \in L^p(0,T;W^{1,q}(\Omega))$ with $p,q>2$,
\begin{eqnarray} \label{eqn:gui}
&&\left \arrowvert \int_0^T\int_{\Omega}g_i(u_i)\partial_t u_i\phi dxdt\right \arrowvert=\left \arrowvert \int_0^T\int_{\Omega}M_i(u_i)\nabla \frac{\mu_i}{g_i(u_i)}\nabla \phi dxdt\right \arrowvert \\ \notag
&\leq& \int_0^T\left ( \left \Arrowvert \sqrt{M_i(u_i)}\nabla \frac{\mu_i}{g_i(u_i)}\right \Arrowvert_{L^2(\Omega)}\left \Arrowvert \sqrt{M_i(u_i)}\right \Arrowvert_{L^{\frac{2q}{q-2}}(\Omega)}\left \Arrowvert \nabla \phi \right \Arrowvert_{L^q(\Omega)}\right)dt\\ \notag
&\leq&\left \Arrowvert M_i(u_i)\right \Arrowvert^{\frac{1}{2}}_{L^{\frac{p}{p-2}}(0,T;L^{\frac{q}{q-2}}(\Omega))} \left \Arrowvert \sqrt{M_i(u_i)}\nabla \frac{\mu_i}{g_i(u_i)}\right \Arrowvert_{L^2(\Omega_T)}\left \Arrowvert \nabla \phi \right \Arrowvert_{L^p(0,T;L^q(\Omega))} \\ \notag
&\leq & C  \left \Arrowvert \phi \right \Arrowvert_{L^p(0,T;W^{1,q}(\Omega))}.
\end{eqnarray}
Let
\begin{equation}
G_i(u_i)=\int_0^{u_i} g_i(a) da. \label{eqn:Gi}
\end{equation}

Thus  \eqref{eqn:gui} yields $\partial_t G_i(u_i)=g_i(u_i)\partial_t u_i \in L^{p'}(0,T; (W^{1,q}(\Omega))')$ with $p'=\frac{p}{p-1}$ and
\begin{equation}
\left \Arrowvert \partial_t G_i(u_i) \right \Arrowvert_{L^{p'}(0,T; (W^{1,q}(\Omega))')} \leq C \text { for all } i. \label{eqn:Git}
\end{equation}
Moreover, by growth assumption on $g$ and estimates on $u_i$, we have
 \begin{equation}
\left \Arrowvert G_i(u_i)\right\Arrowvert_{ L^{\infty}(0,T;W^{1,s}(\Omega))} \leq C. \label{eqn:Gibd}
\end{equation}
 for $1\leq s <2$.
By  \eqref{eqn:uibd}, \eqref{eqn:Mui}, \eqref{eqn:Git}-\eqref{eqn:Gibd} and Remark \ref{rmk:cpt} we can  find a subsequence, not relabeled, a function $u\in L^{\infty}(0,T;H^1(\Omega))$,  a function $\xi \in L^2(\Omega_T)$  and a function $\eta\in L^{\infty}(0,T;W^{1,s}(\Omega)) $ such that as $i \rightarrow \infty$,
\begin{eqnarray}\label{eqn:uiwC}
&&u_i \rightharpoonup u \text{ weakly-* in } L^{\infty}(0,T;H^1(\Omega)), \\ \label{eqn:uip}
&&\sqrt{M_i(u_i)}\nabla \frac{\mu_i}{g_i(u_i)} \rightharpoonup \xi \text{ weakly in } L^2(\Omega_T),\\
&&G_i(u_i)\rightharpoonup \eta\text{ weakly-* in } L^{\infty}(0,T;W^{1,s}(\Omega))\\\label{eqn:GiuiSCp}
&&G_i(u_i) \rightarrow \eta \text{ strongly in } L^{\alpha}(0,T; L^{\beta}(\Omega)) \text{ and a.e. in } \Omega_T,\\ \label{eqn:GuiSCinf}
&&G_i(u_i) \rightarrow \eta \text{ strongly in } C(0,T; L^{\beta}(\Omega)),\\
&& \partial_t G_i(u_i) \rightharpoonup \partial_t \eta \text{ weakly in } L^{p'}(0,T;(W^{1,q}\Omega))').
\end{eqnarray}
where  $1\leq \alpha,\beta <\infty$. By \eqref{eqn:GuiSCinf} and \eqref{eqn:unicov} from Remark \ref{rmk:cptns}, we have
$$
\left \Arrowvert G_i(u_i(x,t+h))-G_i(u_i(x,t))\right\Arrowvert_{C([0,T];L^{\beta}(\Omega))} \rightarrow 0 \text{ uniformly in } i \text{ as } h\rightarrow 0.
$$ Thus given any $\e >0$, there exists $h_{\e}>0$ such that for all $0<h<h_{\e}$ and all $i$,

$$
\left \Arrowvert G_i(u_i(x,t+h))-G_i(u_i(x,t))\right\Arrowvert_{C([0,T];L^{\beta}(\Omega))}^{\beta} <\e.
$$
Given any $\delta>0$, let $I_{\delta}=(1-\delta,1+\delta)\cup (-1-\delta,-1+\delta)$. Consider  the interval having  $u_i(x,t)$ and $u_i(x,t+h)$ as end points.  Denote this interval by $J_i(x,t;h)$. We consider three cases.

\textbf{ Case I: $J_i(x,t;h)\cap I_{\delta}=\varnothing$.}

In this case, $g_i(s)\geq \max\{\theta_i^m, \delta^{2m}\}$ for any $s\in J_i(x,t;h)$ and by \eqref{eqn:Gi}
$$
\left |G_i(u_i(x,t+h))-G_i(u_i(x,t))\right|=\left|\int_{u_i(x,t)}^{u_i(x,t+h)}g_i(s)ds \right|\geq \delta^{2m}|u_i(x,t+h)-u_i(x,t)|.
$$

\textbf{ Case II: $J_i(x,t;h)\cap I_{\delta}\neq \varnothing$ and $|u_i(x,t+h)-u_i(x,t)|\geq 3\delta$ .}

In this case, we have $$|J_i(x,t;h)\cap I_{\delta}^{c}|\geq \frac{1}{3}|J_i(x,t;h)|$$ and
\begin{eqnarray*}
\left |G_i(u_i(x,t+h))-G_i(u_i(x,t))\right|&\geq& \left|\int_{J_i(x,t;h)\cap I_{\delta}^{c}}g_i(s)ds \right|\\
&\geq& \frac{\delta^{2m}}{3}|u_i(x,t+h)-u_i(x,t)|.
\end{eqnarray*}

\textbf{ Case III: $J_i(x,t;h)\cap I_{\delta}\neq \varnothing$ and $|u_i(x,t+h)-u_i(x,t)|< 3\delta$ }

In this case, we have
$$
g_i(s)\leq \max\{(8\delta+16\delta^2)^m,\theta_i^m\} \text{ for any } s \in J_i(x,t;h).
$$Thus
\begin{eqnarray*}
\left |G_i(u_i(x,t+h))-G_i(u_i(x,t))\right|\leq 3\delta\max\{(8\delta+16\delta^2)^m,\theta_i^m\}.
\end{eqnarray*}

Pick $\delta=\e^{\frac{1}{4m\beta}}$ and fix $t$.  Let
$$
\Omega_i^t=\{x\in \Omega: J_i(x,t:h) \text{ satisfies case I or II} \}.
$$ Then
\begin{eqnarray*}
&&\int_{\Omega}\left \arrowvert u_i(x,t+h)-u_i(x,t) \right\arrowvert^{\beta} dx \\
&=&\int_{\Omega_{i}^t}\left \arrowvert u_i(x,t+h)-u_i(x,t) \right\arrowvert^{\beta} dx+\int_{\Omega\backslash\Omega_{i}^t}\left \arrowvert u_i(x,t+h)-u_i(x,t) \right\arrowvert^{\beta} dx\\
&\leq& 3^{\beta}\e^{-\frac{1}{2}}\int_{\Omega_{i}^t}\left \arrowvert G_i(u_i(x,t+h))-G_i(u_i(x,t)) \right\arrowvert^{\beta} dx+\int_{\Omega\backslash\Omega_{i}^t}\left \arrowvert u_i(x,t+h)-u_i(x,t) \right\arrowvert^{\beta} dx\\
&\leq& 3^{\beta}\e^{\frac{1}{2}}+C\e^{\frac{1}{4m}}
\end{eqnarray*}
Taking maximum on the left side, we have for all $i$, any $h<h_{\e}$,
$$
\left \Arrowvert (u_i(x,t+h)-u_i(x,t)\right\Arrowvert_{C([0,T];L^{\beta}(\Omega))}^{\beta} \leq 3^{\beta}\e^{\frac{1}{2}}+C\e^{\frac{1}{4m}}.
$$

Thus
$$
\left \Arrowvert u_i(x,t+h)-u_i(x,t)\right\Arrowvert_{C([0,T];L^{\beta}(\Omega))}^{\beta} \rightarrow 0 \text{ uniformly as } h\rightarrow 0.
$$

In addition, for any $0<t_1<t_2<T$, \eqref{eqn:uibd} implies that  for $1\leq \beta <\infty$, we have
$$
\int_{t_1}^{t_2} u_i(x,t)dt \text{ is relatively compact in }L^{\beta}(\Omega).
$$
Therefore we conclude from Remark \ref{rmk:cptns} that
\begin{equation}
u_i\rightarrow u(x,t) \text{ strongly in  }C([0,T];L^{\beta}(\Omega)) \text{ for } 1\leq \beta <\infty.\label{eqn:uistronginf}
\end{equation}
Similarly. we can prove
\begin{equation}
u_i\rightarrow u(x,t) \text{ strongly in  }L^{\alpha}(0,T;L^{\beta}(\Omega)) \text{ for } 1\leq \alpha,\beta <\infty \text{ and a.e. in }\Omega_T.\label{eqn:uistronglp}
\end{equation}
Growth condition on $M(u)$ and \eqref{eqn:uistronginf}, \eqref{eqn:uistronglp} yield
\begin{eqnarray*}
&&M_i(u_i)\rightarrow M(u) \text{ strongly in }C([0,T];L^{\beta}(\Omega)) \text{ for } 1\leq \beta <\infty, \\
&&M_i(u_i)\rightarrow M(u) \text{ strongly in }L^{\alpha}(0,T;L^{\beta}(\Omega))  \text{ for } 1\leq \alpha,\beta <\infty,\\&&\sqrt{M(_i(u_i))}\rightarrow \sqrt{M(u)} \text{ strongly in  }C([0,T];L^{\gamma}(\Omega)) \text{ for } 1\leq \gamma <\infty.
\end{eqnarray*}
Hence $G_i(u_i)$ converges to $G(u)$ a.e. in $\Omega_T$ and $\eta=G(u)$.
Passing to the limit in \eqref{eqn:ui2}, we have
\begin{equation}
\int_0^T\left<g(u)\partial_tu, \phi\right>_{((W^{1,q}(\Omega))',W^{1,q}(\Omega))} dt=-\int_0^T\int_\Omega \sqrt{M(u)}\xi\cdot \nabla \phi dxdt \label{eqn:ueq}
\end{equation}
for any $\phi \in L^p(0,T;W^{1,q}(\Omega))$ with $p,q>2$.
\begin{remark} \label{rmk:cptns}
(Compactness in $L^p(0,T;B)$ Theorem 1 in \cite{Sim86})
Assume $B$ is a Banach space and $F\subset L^p(0,T;B)$. $F$ is relatively compact in  $L^p(0,T;B)$ for $1\leq p<\infty$, or in $C([0,T],B)$ for $p=\infty$ if and only if
\begin{equation}
\left\{\int_{t_1}^{t_2}f(t)dt: f\in F\right \} \text{ is relatively compact in }B, \forall   0<t_1<t_2<T, \label{eqn:unibd}
\end{equation}
\begin{equation}
\left\Arrowvert\tau_hf-f\right \Arrowvert_{L^p(0,T;B)}\rightarrow 0 \text{ as } h\rightarrow 0 \text{ uniformly for } f\in F. \label{eqn:unicov}
\end{equation}
Here $\tau_h f(t)=f(t+h)$ for $h>0$ is defined on $[-h,T-h]$.
\end{remark}

\subsection{Weak convergence of  $\nabla \frac{\mu_i}{g_i(u_i)}$}
We now look for relation between $\xi$ and $u$. Following the  idea in \cite{DaiDu16}, we  decompose $\Omega_T$ as follows. Let  $\delta_j$ be a positive sequence monotonically decreasing to $0$. By \eqref{eqn:uip} and Egorov's theorem, for every $\delta_j>0$, there exists $B_j \subset \Omega_T$ satisfying $|\Omega_t\backslash B_j|< \delta_j$ such that
\begin{equation}
u_i\rightarrow u \text{ uniformly in } B_j. \label{eqn:unifui}
\end{equation}
We can pick
\begin{equation}
B_1\subset B_2 \subset \cdots \subset B_j\subset B_{j+1} \subset\cdots \subset \Omega_T. \label{eqn:Bj}
\end{equation}
Define
$$
P_j:=\{(x,t)\in \Omega_T:|1-u^2| > \delta_j\}.
$$
Then
\begin{equation}
P_1\subset P_2 \subset \cdots \subset P_j \subset P_{j+1}\subset \cdots \subset\Omega_T. \label{eqn:Pj}
\end{equation}
Let  $B=\cup _{j=1}^{\infty}B_j$ and $P=\cup_{j=1}^{\infty}P_j$.  Then $|\Omega_T\backslash B|=0$ and each $B_j$ can be split into two parts:
\begin{eqnarray*}
&&D_j=B_j\cap P_j, \text{ where } |1-u^2|>\delta_j, \text{ and } u_i\rightarrow u \text{ uniformly},\\
&&\tilde D_j=B_j\backslash P_j, \text{ where } |1-u^2|\leq \delta_j, \text{ and } u_i \rightarrow u \text { uniformly }.
\end{eqnarray*}
\eqref{eqn:Bj} and \eqref{eqn:Pj} imply
\begin{equation}
D_1\subset D_2 \subset \cdots \subset D_j \subset D_{j+1}\subset \cdots \subset D:=B\cap P. \label{eqn:Dj}
\end{equation}
For any $\Psi \in L^p(0,T;L^q(\Omega,\R^n))$ with $ p,q>2$, we have
\begin{eqnarray}\notag
&&\int_{\Omega_T}M_i(u_i)\nabla \frac{\mu_i}{g_i(u_i)}\cdot \Psi dxdt\\ \notag
&=&\int_{\Omega_T\backslash B_j}M_i(u_i)\nabla \frac{\mu_i}{g_i(u_i)}\cdot \Psi dxdt+\int_{D_j}M_i(u_i)\nabla \frac{\mu_i}{g_i(u_i)}\cdot \Psi dxdt\\ \label{eqn:rhsi}
&&+\int_{\tilde D_j}M_i(u_i)\nabla \frac{\mu_i}{g_i(u_i)}\cdot \Psi dxdt
\end{eqnarray}
The left hand side of \eqref{eqn:rhsi} converges to $\int_{\Omega_T}\sqrt{M(u)}\xi\cdot \Psi dxdt$. We analyze the three terms on the right hand side separately.
 To estimate the first term on the right hand side of \eqref{eqn:rhsi}, noticing $|\Omega_T\backslash B_j| \rightarrow 0$ and
\begin{eqnarray*}
\lim_{i\rightarrow \infty} \int_{\Omega_T\backslash B_j}M_i(u_i)\nabla \frac{\mu_i}{g_i(u_i)}\cdot \Psi dxdt=\ \int_{\Omega_T\backslash B_j}\sqrt{M(u)}\xi\cdot \Psi dxdt,
\end{eqnarray*}
we have
\begin{eqnarray*}
\lim_{j\rightarrow \infty}\lim_{i\rightarrow \infty}\int_{\Omega_T\backslash B_j}M_i(u_i)\nabla \frac{\mu_i}{g_i(u_i)}\cdot \Psi dxdt=0.
\end{eqnarray*}
By uniform convergence of $u_i$ to $u$ in $B_j$, we introduce subsequence $u_{j,k}$ such that $u_{j,k}\rightarrow u$ uniformly in $B_j$ and there exists $N_j $ such that for all $k\geq N_j$,
\begin{equation}
|1-u^2_{j,k}|>\frac{\delta_j}{2} \text{ in } D_j,    \hspace{0.2 in} |1-u^2_{j,k}|\leq 2\delta_j \text{ in } \tilde D_j.
\end{equation}
Thus the third term on the right hand side of \eqref{eqn:rhsi} can be estimated by
\begin{eqnarray*}
&&\lim_{j\rightarrow \infty}\lim_{k\rightarrow \infty}\left \arrowvert \int_{\tilde D_j}M_{j,k}(u_{j,k})\nabla \frac{\mu_{j,k}}{g_{j,k}(u_{j,k})} \cdot \Psi dxdt \right \arrowvert \\
&\leq&\lim_{j\rightarrow \infty}\lim_{k\rightarrow \infty}\left \{\left(\sup_{\tilde D_j}\sqrt{M_{j,k}(u_{j,k})}\right)\left \Arrowvert \Psi \right \Arrowvert_{L^2(\tilde D_j)}\left \Arrowvert \sqrt{M_{j,k}(u_{j,k})}\nabla \frac{\mu_{j,k}}{g_{j,k}(u_{j,k})}\right \Arrowvert_{L^2(\tilde D_j)}\right \} \\
&\leq&\left(\sup_{\tilde D_j}\sqrt{M_{j,k}(u_{j,k})}\right)|\Omega|^{\frac{q-2}{2q}}\left \Arrowvert \Psi \right \Arrowvert_{L^2(0,T;L^q(\Omega)}\left \Arrowvert \sqrt{M_{j,k}(u_{j,k})}\nabla \frac{\mu_{j,k}}{g_{j,k}(u_{j,k})}\right \Arrowvert_{L^2(\tilde D_j)}\\
&&\leq C\lim_{j\rightarrow \infty}\lim_{k\rightarrow \infty}\max\left \{(2\delta_j)^{m/2}, \theta_{j,k}^{m/2}\right\}\\
&=&0.
\end{eqnarray*}
For the second  term, we see that
\begin{eqnarray*}
&&\left(\frac{\delta_j}{2}\right)^m\int_{D_j}|\nabla \frac{\mu_{j,k}}{g_{j,k}(u_{j,k})}|^2 dxdt \\
&\leq&\int_{D_j}M_{j,k}(u_{j,k})|\nabla \frac{\mu_{j,k}}{g_{j,k}(u_{j,k})}|^2 dxdt\\
&\leq& \int_{\Omega_T}M_{j,k}(u_{j,k})|\nabla \frac{\mu_{j,k}}{g_{j,k}(u_{j,k})}|^2 dxdt \leq C.
\end{eqnarray*}
Therefore $\nabla \frac{\mu_{j,k}}{g_{j,k}(u_{j,k})}$ is bounded in $L^2(D_j)$ and we can extract a further subsequence, not relabeled, which converges weakly to some $\xi_j \in L^2(D_j)$.
Since $D_j$ is an increasig sequence of sets with $\lim_{j\rightarrow \infty} D_j=D$, we have $\xi_j=\xi_{j-1}$ a.e. in $D_{j-1}$. By setting $\xi_j=0$ outside $D_j$, we can extend $\xi_j $ to a $L^2$  function $\tilde \xi_j$ defined in $D$. Therefore for a.e. $x\in D$, there exists a limit of $\tilde \xi_j(x)$ as $j \rightarrow \infty$. Let $\xi(x)=\lim_{j\rightarrow \infty}\tilde \xi_j(x)$, we see that $\xi(x)=\xi_j(x)$ for a.e $x \in D_j$ and for all $j$.

By a standard diagonal  argument, we can extract a subsequnce such that
\begin{equation}
\nabla \frac{\mu_{k,N_k}}{g_{k,N_k}(u_{k,N_k})}\rightharpoonup \zeta \text{ weakly in } L^2(D_j) \text{ for all } j. \label{eqn:mugwC}
\end{equation}

By strong convergence of $\sqrt{M_i(u_i)}$ to $\sqrt{M(u)}
$ in $C([0,T];L^{\beta}(\Omega))$ for $1\leq \beta<\infty$, we obtain
$$
\chi_{D_j}\sqrt{M_{k,N_k}(u_{k,N_k})}\nabla \frac{\mu_{k,N_k}}{g_{k,N_k}(u_{k,N_k})}\rightharpoonup \chi_{D_j}\sqrt{M(u)}\zeta
$$
weakly in $L^2(0,T;L^q(\Omega))$ for $1\leq q <2$ and all $j$.  Recall $\sqrt{M_i(u_i)}\nabla \frac{\mu_i}{g_i(u_i)} \rightarrow \xi$ weakly in $L^2(\Omega_T)$, we have $\xi=\sqrt{M(u)}\zeta$ in $D_j$ for all $j$. Hence $\xi=\sqrt{M(u)}\zeta$ in $D$ and consequently
$$
\chi_DM_{k,N_k}(u_{k,N_k})\nabla \frac{\mu_{k,N_k}}{g_{k,N_k}(u_{k,N_k})}\rightharpoonup \chi_{D}M(u)\zeta
$$
weakly in $L^2(0,T;L^q(\Omega))$ for $1\leq q <2$.

Replacing $u_i$ by subsequence $u_{k,N_k}$ in \eqref{eqn:rhsi} and letting $k \rightarrow \infty$ then $j\rightarrow \infty$, we have
\begin{eqnarray}\notag
\int_{\Omega_T}\sqrt{M(u)}\xi \cdot \Psi dxdt&=&\lim_{j\rightarrow \infty}\int_{D_j}M(u)\zeta \cdot \Psi dxdt \\ \label{eqn:limit}
&=&\int_D M(u)\zeta\cdot \Psi dxdt.
\end{eqnarray}
It follows from \eqref{eqn:ueq} and \eqref{eqn:limit} that
\begin{equation}
\int_0^T\left<g(u)\partial_t u,\phi\right>_{((W^{1,q}(\Omega))',W^{1,q}(\Omega))} dt=-\int_{D}M(u)\zeta\cdot \nabla \phi dxdt
\end{equation}
for all $\phi \in L^p(0,T;W^{1,q}(\Omega))$ where $p,q>2$.

\subsection{Relation between $\zeta$ and $u$}

The desired relation between $\zeta$ and $u$ is
\begin{eqnarray}\notag
\zeta&=&\frac{1}{g}\nabla \mu-\mu \frac{g'(u)}{g^2(u)}\nabla u \\ \label{eqn:zeta}
&=&\frac{1}{g}\left(-\nabla \Delta u+q''(u)\nabla u\right)-\frac{g'(u)}{g^2(u)}\nabla u \left (-\Delta u+q'(u)\right).
\end{eqnarray}
Given the known regularity $u\in L^{\infty}(0,T;H^1(\Omega))$ and degeneracy of $g(u)$,  the right hand side of \eqref{eqn:zeta} might not be defined as a function.   We can, however,
under the additional assumption $u\in L^2(0,T;H^2(\Omega))$ and suitable assumptions on integrability of $\nabla \Delta u$, find an explicit expression of $\zeta$ in terms of \eqref{eqn:zeta} in suitable subset of $\Omega_T$.

{\it{Claim I: If $u \in L^2(0,T:H^2(\Omega))$  and for some $j$, the interior of $D_j$, denoted by $(D_j)^{\circ}$, is not empty, then
$$
\nabla \Delta u \in L^1((D_j)^{\circ})
$$
and
$$
\zeta=\frac{-\nabla \Delta u+q''(u)\nabla u}{g(u)}-\frac{g'(u)}{g^2(u)}\left(-\Delta u+q'(u) \right)\nabla u \text{ a.e. in } (D_j)^{\circ}.
$$}}

{\it{Proof of the claim I.}} Since $u \in L^2(0,T; H^2(\Omega))$, we can have a subsequence, not relabeled such that,   $u_{k,N_k}$ converges weakly to $ u$ in $L^2(0,T; H^2(\Omega))$. Since
\begin{equation}
\mu_{k,N_k}=-\Delta u_{k,N_k}+q'(u_{k,N_k})\text{ in } \Omega_T, \label{eqn:muik}
\end{equation}
The right hand side of \eqref{eqn:muik} weakly converges to $-\Delta u+q'(u)$ in $L^2(\Omega_T)$. Hence
$$
\mu_{k,N_k} \rightharpoonup \mu=-\Delta u+q'(u) \text{ weakly in }  L^2(\Omega_T).
$$
On the other hand, using  $u_{k,N_k}$ and $u$ as test functions  in  \eqref{eqn:mutheta} yield
\begin{eqnarray*}
&&\int_{\Omega_T}\mu_{k,N_k}u_{k,N_k}dxdt=\int_{\Omega_T}\left (\left \arrowvert \nabla u_{k,N_k}\right\arrowvert^2+q'(u_{k,N_k})u_{k,N_k}\right)dxdt\\
&&\int_{\Omega_T}\mu_{k,N_k}udxdt=\int_{\Omega_T}\left (\nabla u_{k,N_k}\cdot \nabla u+q'(u_{k,N_k})u \right)dxdt.
\end{eqnarray*}
Passing to the limit, by \eqref{eqn:uistronglp}, growth assumptions on $q'$ and \eqref{eqn:muik}, we have
$$
\lim_{k\rightarrow \infty} \int_{\Omega_T}\left \arrowvert \nabla u_{k,N_k}\right\arrowvert^2=\int_{\Omega_T}\left \arrowvert \nabla u\right\arrowvert^2.
$$
Therefore $$\nabla u_{k,N_k} \rightarrow \nabla u \text{ strongly in } L^2(\Omega_T).$$
Since $u_{k,N_k} \in L^2(0,T;W^{3,s}(\Omega))$,  we can differentiate \eqref{eqn:muik} and get
\begin{equation}
\nabla \mu_{k,N_k}=-\nabla \Delta u_{k,N_k}+q''(u_{k,N_k})\nabla u_{k,N_k}, \label{eqn:gmuk}
\end{equation}
and
\begin{equation}
\nabla \frac{\mu_{k,N_k}}{g_{k,N_k}(u_{k,N_k})}=\frac{1}{g_{k,N_k}(u_{k,N_k})}\nabla \mu_{k,N_k}-\mu_{k,N_k}\frac{g'_{k,N_k}(u_{k,N_k})}{g^2_{k,N_k}(u_{k,N_k})}\nabla u_{k,N_k} \label{eqn:muikuik}
\end{equation}
on $D_j^{\circ}$.
Thus
\begin{equation}
\nabla \mu_{k,N_k}=g_{k,N_k}(u_{k,N_k})\nabla  \frac{\mu_{k,N_k}}{g_{k,N_k}(u_{k,N_k})}+ \frac{\mu_{k,N_k}}{g_{k,N_k}(u_{k,N_k})}g'_{k,N_k}(u_{k,N_k})\nabla u_{k,N_k}.\label{eqn:gmuik}
\end{equation}
Since
\begin{eqnarray*}
&& g_{k,N_k}(u_{k,N_k})\rightarrow g(u) \text{ uniformly in } D_j^{\circ},\\
&& \frac{g'_{k,N_k}(u_{k,N_k})}{g_{k,N_k}(u_{k,N_k})} \rightarrow \frac{g'(u)}{g(u)} \text{ uniformly in } D_j^{\circ},\\
&& \nabla  \frac{\mu_{k,N_k}}{g_{k,N_k}(u_{k,N_k})} \rightharpoonup \zeta  \text{ weakly in } L^2(D_j^{\circ}),\\
&&\mu_{k,N_k}\rightharpoonup \mu \text{ weakly in } L^2(\Omega_T),\\
&&\nabla u_{k,N_k}\rightarrow \nabla u \text{ strongly in } L^2(\Omega_T),
\end{eqnarray*}
we have, for any $\phi \in L^{\infty}(D_j^{\circ})$,
\begin{eqnarray*}
&&\int_{D_j^{\circ}} \phi \left(g_{k,N_k}(u_{k,N_k})\nabla  \frac{\mu_{k,N_k}}{g_{k,N_k}(u_{k,N_k})}+ \frac{\mu_{k,N_k}}{g_{k,N_k}(u_{k,N_k})}g'_{k,N_k}(u_{k,N_k})\nabla u_{k,N_k}\right) dxdt \\
&&\rightarrow \int_{D_j^{\circ}}\phi \left (g(u)\zeta+\frac{g'(u)}{g(u)}\mu \nabla u\right) dxdt,
\end{eqnarray*}
i.e.
$$
\nabla \mu_{k,N_k}\rightharpoonup \eta\coloneq g(u)\zeta+\frac{g'(u)}{g(u)}\mu \nabla u\ \text{ weakly in } L^1(D_j^{\circ}).
$$
Passing to the limit in \eqref{eqn:gmuk},  we obtain, in  the sense of distribution, that
$$
\eta=-\nabla \Delta u+q''(u)\nabla u.
$$

Since $q''(u)\nabla u\in L^2(\Omega_T)$, we have  $-\nabla \Delta u \in L^1(D_j^{\circ})$, hence
\begin{equation}
\eta=-\nabla \Delta u+q''(u)\nabla u \text{ a.e. in } D_j^{\circ}
\end{equation}

Since $\frac{1}{g_{k,N_k}(u_{k,N_k})} \rightarrow \frac{1}{g(u)}$ uniformly  in $D_j$, we have
$$
\frac{1}{g_{k,N_k}}\nabla \mu_{k,N_k} \rightharpoonup \frac{1}{g(u)}\eta \text{ weakly in } L^1(D_j^{\circ}).
$$

Since $\frac{g'_{k,N_k}(u_{k,N_k})}{g^2_{k,N_k}(u_{k,N_k})} \rightarrow \frac{g'(u)}{g^2(u)}$ uniformly  in $D_j$, we have
$$
\frac{g'_{k,N_k}(u_{k,N_k})}{g^2_{k,N_k}(u_{k,N_k})}\mu_{k,N_k}\nabla u_{k,N_k} \rightharpoonup \frac{g'(u)}{g^2(u)}\mu \nabla u \text{ weakly in } L^1(D_j^{\circ}).
$$

Passing to the limit in \eqref{eqn:muikuik}, we have
\begin{eqnarray*}
\zeta&=&\frac{1}{g(u)}\eta-\mu \frac{g'(u)}{g^2(u)}\nabla u=\frac{-\nabla \Delta u+q''(u)\nabla u}{g(u)}-\frac{g'(u)}{g^2(u)}\left(-\Delta u+q'(u)\right)\nabla u
\end{eqnarray*}
on $(D_j)^{\circ}$.
Noticing the value of $\zeta$  on $\Omega_T\backslash D$ doesn't matter since it does not appear on the right hand side of \eqref{eqn:limit}.

{\it{Claim II: For any open set $U\in \Omega_T$ in which  $\nabla \Delta u \in L^p(U)$ for some $p>1$ and $g(u)>0$, we have
\begin{equation}
\zeta=\frac{-\nabla \Delta u+q''(u)\nabla u}{g(u)}-\frac{g'(u)}{g^2(u)}\left(-\Delta u+q'(u) \right)\nabla u . \label{eqn:zetaexp}
\end{equation} in $U$.}}

To prove this, since
\begin{equation}
\nabla \mu_{k,N_k}=-\nabla \Delta u_{k,N_k}+q''(u_{k,N_k})\nabla u_{k,N_k} \text{ in } \Omega_T \label{eqn:muknk}
\end{equation}
and
\begin{equation}
\nabla \frac{\mu_{k,N_k}}{g_{k,N_k}(u_{k,N_k})} =\frac{1}{g_{k,N_k}(u_{k,N_k})}\nabla \mu_{k,N_k}+\mu_{k,N_k}\cdot \nabla \frac{1}{g_{k,N_k}(u_{k,N_k})} \text{ on } D_j\label{eqn:mugk}.
\end{equation}

The right hand side of \eqref{eqn:muknk} converges weakly to $-\nabla \Delta u+q''(u)\nabla u$ in $L^q(U)$ for $q=\min\{p,2\}>1$. Hence
$$
\nabla \mu_{k,N_k} \rightharpoonup \eta=-\nabla \Delta u+q''(u)\nabla u \text{ weakly in } L^q(U).
$$

The right hand side of \eqref{eqn:mugk} converges weakly to
$$\frac{\eta}{g(u)} -\frac{g'(u)}{g^2(u)}\mu\cdot \nabla u$$
in $L^1(U\cap D_j)$ for each $j$ and therefore $$\zeta=\frac{-\nabla \Delta u+q''(u)\nabla u}{g(u)}-\frac{g'(u)}{g^2(u)}\left(-\Delta u+q'(u) \right)\nabla u $$ in $ U\cap D$.  The definition of $\zeta$ can be extended to $U\backslash D$ by our integrability assumption on $u$. Define
$$
\tilde \Omega_T=\{U\subset \Omega_T: g(u)>0 \text{ on } U \text{ and }\nabla \Delta u \in L^p(U) \text{ for some } p>1 \text{ depending on } U\}.
$$
Then $\tilde \Omega_T$ is open and $\zeta$ is defined by  \eqref{eqn:zetaexp} on  $\tilde \Omega_T$. Since $|\Omega_T\backslash B|=0$ ,  $M(u)=0$ on $\Omega_T\backslash P$ and
$$
\Omega_T\backslash \{D\cup \tilde \Omega_T\}\subset \{\Omega_T\backslash B\} \cup \{\Omega_T\backslash P\},
$$
we can take the value  of $\zeta$ to be zero outside $D\cup \Omega_T$, and it won't affect the integral on the right side of \eqref{eqn:dg}.

Lastly the energy inequality \eqref{eqn:ineq}  follows by taking limit in  the energy inequality for $u_{k,N_k}$.

\begin{remark}
 In  Cahn-Hilliard case,  there is convergence of $\nabla \mu_k$ on $L^2(D_j)$, and  relation between $\xi$ and $u$ can be derived directly. Here we only have convergence of $\nabla \frac{\mu_k}{g_k(u_k)}$  on $L^2(D_j)$. In order to obtain convergence of $\nabla \mu_k$, we need convergence  $\mu_k$ on $L^p(\Omega_T)$ for suitable $p$, this is where we used the additional assumption $u\in L^2(0,T;H^2(\Omega))$.

\end{remark}
\section{A Modified phase field model for self-climb of prismatic dislocation loops} \label{sec:pf}

Dislocations are line defects in crystals \cite{Hirth-Lothe,Xiang2006}.
A phase field model \cite{NXY20} was derived based on the pipe diffusion model for self-climb of prismatic dislocation loops \cite{Niu2019,Niu2017} that describes the conservative climb of dislocation loops observed in experiments of irradiated materials \cite{Kroupa,Hirth-Lothe,Dudarev2013}.
In this section, we study the wellposedness of the following modified phase field model for self-climb of prismatic dislocation loops:
\begin{eqnarray}
\label {eqn:rrv} g(u)\partial_t u&=&\nabla \cdot (M(u)\nabla \frac{\mu}{g(u)}) \text{ for } x \in \Omega \subset \mathbb{R}^2, t\in [0,\infty) \\
\label{eqn:chemrrv} \mu&=&-\Delta u+\frac{1}{\e^2}q'(u)+\frac{1}{\e}f_{cl}
\end{eqnarray}

Where $M(u)=M_0g(u)$, $ g(u)=|1-u^2|^m$ for $2\leq m<\infty$, $q(u)$ satisfy same assumptions \eqref{eqn:qu}-\eqref{eqn:q'u} as those for Eqs.~\eqref{eqn:ch}-\eqref{eqn:chem}. Here $f_{cl}$ is the total climb force with
$$
f_{cl}=f_{cl}^d+f_{cl}^{app}
$$
where  $f_{cl}^{app}$ is the applied climb force, and
\begin{equation}
f_{cl}^d(x,y,u)=\frac{G b^2}{4\pi(1-\nu)}\int_{\Omega}\left(\frac{x-\ov x}{R^3}u_{\ov x}+\frac{y-\ov{y}}{R^3}u_{\ov{y}}\right)d\ov x d\ov y \label{equ: clfrrv}
\end{equation}
represents the climb force generated by all the dislocations. Here $\Omega \subset \R^2$ is a bounded domain, $G$ is the shear modulus, $\nu$ is the Poisson ratio, and $R=\sqrt{(x-\ov x)^2+(y-\ov y)^2}$. In this model, we assume that the prismatic dislocation loops lie and evlove by self-climb in the $xy$ plane and all dislocation loops have the same Burgers vector $\mathbf{b}=(0,0,b)$.

The chemical potential $\mu$  comes from variations of the classical Cahn-Hilliard energy and the elastic energy due to dislocations, i.e.
\begin{equation}
\mu=\frac{\delta E_{CH}}{\delta u}+\frac{\delta E_{el}}{\delta u},
\end{equation}
where
\begin{eqnarray}
&&E_{CH}(u)=\int_{\Omega}\left (\frac{1}{2}|\nabla u|^2+q(u)\right) dx, \\
&& E_{el}=\int_{\Omega}\left ( \frac{1}{2}u f^{d}_{cl}+u f_{cl}^{app}\right ) dx
\end{eqnarray}
are classical Cahn-Hilliard energy and elastic energy, respectively.
Under periodic boundary conditions,  the climb force generated by the dislocations can be expressed as
\begin{eqnarray}
f^{d}_{cl}(x,y,u)=\frac{G b^2}{2(1-\nu)}(-\Delta)^{\frac{1}{2}}u.
\end{eqnarray}
Here $(-\Delta)^{s}u$ is a fractional operator defined by $$\mathcal{F}((-\Delta)^s f)=(\xi_1^2+\xi_2^2)^{\frac{s}{2}}\mathcal{F}(f)(\xi)$$ for $\xi \in \mathbb{Z}^2$. In the analysis below, without loss of generality, we set the coefficient of the climb force $\frac{G b^2}{2(1-\nu)}=1$.

System  \eqref{eqn:rrv}-\eqref{eqn:chemrrv} is a modified version of the phase field model introduced in \cite{NXY20}, which does not have the $g(u)$ term on the left side of \eqref{eqn:rrv}. Putting an extra factor $h=H_0g$ in front of the nonlocal climb force $f^{d}_{cl}$, the asymptotic analysis in \cite{NXY20} showed that the proposed phase field model yields accurate dislocation self-climb velocity in  the sharp interface limit. Moreover, numerical simulations in \cite{NXY20} showed excellent agreement with experimental observations and discrete dislocation dynamics simulation results. Now we prove the wellposedness of the  modified model \eqref{eqn:rrv}-\eqref{eqn:chemrrv}. There is an extra nonlocal term $f^{d}_{cl}$ in this model compared with the model considered in previous sections.


Define
$$
u^N(x,t)=\sum_{j=1}^Nc_j^N(t)\phi_j(x), \hspace{0.5 in} \mu^N(x,t)=\sum_{j=1}^Nd^N_{j}(t)\phi_j(x),
$$
where $\{c_j^N,d^N_j\}$ satisfy
\begin{eqnarray}\label{eqn:uNrrv}
\int_{\Omega}\partial_t u^N\phi_j dx&=&-\int_{\Omega}M_{\theta}(u^N)\nabla \frac{\mu^N}{g_{\theta}(u^N)}\cdot \nabla \frac{\phi_j}{g_{\theta}(u^N)}dx,\\ \label{eqn:muNrrv}
\int_{\Omega}\mu^N\phi_j dx&=&\int_{\Omega}\left(\nabla u^N\cdot \nabla \phi_j+q'(u^N)\phi_j+\phi_j(-\Delta)^{\frac{1}{2}} u^N \right)dx,\\ \label{eqn:uN0rrv}
u^N(x,0)&=&\sum_{j=1}^N\left(\int_{\Omega} u_0\phi_j dx\right)\phi_j(x).
\end{eqnarray}
\eqref{eqn:uNrrv}-\eqref{eqn:uN0rrv} is an initial value problem for a system of ordinary equations for $\{c_j^N(t)\}$. Since right hand side of \eqref{eqn:uNrrv} is continuous in $c_j^N$, the system has a local solution.

Define energy functional
$$
F(u)=\int_{\Omega}\left\{\frac{1}{2}|\nabla u|^2 +q(u)+|(-\Delta)^{\frac{1}{4}}u|^2\right\} dx.
$$
Direct calculation yields
$$
\frac{d}{dt}F(u^N(x,t))=-\int_{\Omega}M_{\theta}(u^N)\left|\nabla \frac{\mu^N}{g_{\theta}(u^N)}\right|^2 dx,
$$
integration over $t$ gives the following energy identity for any $t>0$
\begin{eqnarray}\notag
&&\int_{\Omega}\left(\frac{1}{2}|\nabla u^N(x,t)|^2+q(u^N(x,t))+u^N(-\Delta)^{\frac{1}{2}}u^N\right)dx\\ \notag
&&+\int_0^t\int_{\Omega}M_{\theta}(u^N(x,\tau))\left|\nabla \frac{\mu^N(x,\tau)}{g_{\theta}(u^N(x,\tau))}\right|^2 dxd\tau\\ \label{eqn:energyidrrv}
&=&\int_{\Omega}\left(\frac{1}{2}|\nabla u^N(x,0)|^2+q(u^N(x,0))+u^N(x,0)(-\Delta)^{\frac{1}{2}}u^N(x,0)\right)dx\\ \notag
&\leq&\int_{\Omega}\left \Arrowvert \nabla u_0\right \Arrowvert_{L^2(\Omega)}^2+C\left(\left\Arrowvert u_0\right\Arrowvert^{r+1}_{H^1{\Omega}}+|\Omega|\right)+\frac{1}{2}\left\Arrowvert u_0\right\Arrowvert_{L^2(\Omega)}^2\leq C <\infty,
\end{eqnarray}
where  $C$ represents a generic constant possibly depending only on  $T$, $\Omega$,  $u_0$ but not on $\theta$. Since $\Omega$ is bounded region, by  growth assumption assumption  \eqref{eqn:qu} and Poincare's inequality, the energy identity \eqref{eqn:energyidrrv} implies
$
u^N\in L^{\infty}(0,T;H^1(\Omega))
$
with
\begin{equation}
\left\Arrowvert u^N\right\Arrowvert _{L^{\infty}(0,T;H^1(\Omega))}\leq C \text{ for all } N, \label{eqn:uNbdrrv}
\end{equation}
and
\begin{equation}
\left \Arrowvert \sqrt{M_{\theta}(u^N)}\nabla \frac{\mu^N}{g_{\theta}(u^N)}\right \Arrowvert_{L^2(\Omega_T)}\leq C \text{ for all } N.  \label{eqn:muNbdrrv}
\end{equation}

Repeat the argument in Section \ref{sec:ndg} and Section \ref{sec:dg}, replacing  energy functional $F(u)$ by $E(u)$ when necessary, we can prove the following existence theorem for \eqref{eqn:rrv}-\eqref{eqn:chemrrv} with  nondegenerate and degenerate mobilities respectively.

\begin{theorem}\label{thm:ndgpf}
Let  $M_{\theta},  g_{\theta}$ be defined by \eqref{eqn:Mtheta} and \eqref{eqn:gtheta}, under the assumptions \eqref{eqn:qu}-\eqref{eqn:q''u}, for any  $u_0 \in H^1(\Omega)$ and  any $T>0$, there exists a function $u_{\theta}$ such that
\begin{itemize}
\item[a)] $u_{\theta}\in L^{\infty}(0,T;H^1(\Omega))\cap C([0,T];L^p(\Omega))\cap L^2(0,T;W^{3,s}(\Omega))$, where $1\leq p <\infty$, $1\leq s <2$,
\item[b)]$\partial_tu_{\theta}\in L^2(0,T;(W^{1,q}(\Omega))')$ for $q>2$,
\item[c)]$u_{\theta}(x,0)=u_0(x) $ for all $x \in \Omega$,
\end{itemize}
which satisfies \eqref{eqn:u}-\eqref{eqn:mu} in the following weak sense
\begin{eqnarray}\notag
&&\int_0^T<\partial_tu_{\theta},\phi>_{(W^{1,q}(\Omega))',W^{1,q}(\Omega)} dt\\
&&=-\int_0^T\int_{\Omega}M_{\theta}(u_{\theta})\nabla \frac{-\Delta u_{\theta}+q'(u_{\theta})+(-\Delta)^{\frac{1}{2}}u_\theta}{g_{\theta}(u_{\theta})}\cdot \nabla \frac{\phi}{g_{\theta}(u_{\theta})}dxdt
\label{eqn:thetandgrrv}
\end{eqnarray}
for all $\phi\in L^2(0,T;W^{1,q}(\Omega))$ with $q>2$. In addition, the following energy inequality holds for all $t>0$.
\begin{eqnarray}
\notag
&&\int_{\Omega}\left (\frac{1}{2}|\nabla u_{\theta}(x,t)|^2+q(u_{\theta}(x,t))+u_\theta(x,t)(-\Delta)^{\frac{1}{2}}u_\theta\right) dx\\
\notag &&+\int_0^t\int_{\Omega}M_{\theta}(u_{\theta}(x,\tau)\left|\nabla \frac{-\Delta u_{\theta}(x,\tau)+q'(u_{\theta}(x,\tau)+(-\Delta)^{\frac{1}{2}}u_\theta}{g_{\theta}(u_{\theta}(x,\tau))}\right|^2dxd\tau \\
\label{eqn:engineqrrv}&\leq&\int_{\Omega}\left(\frac{1}{2}|\nabla u_0(x)|^2+q(u_{0}(x))+u_0(x)(-\Delta)^{\frac{1}{2}}u_0\right) dx.
\end{eqnarray}
\end{theorem}

\begin{theorem}\label{thm:dgpf}
For any $u_0\in H^1(\Omega)$ and $T>0$, there exists a function $u:\Omega_T=\Omega\times [0,T]\rightarrow \R$ satisfying
\begin{itemize}
\item[i)] $u\in L^{\infty}(0,T;H^1(\Omega))\cap C([0,T];L^s(\Omega))$, where $1\leq s <\infty$,
\item[]
\item[ii)] $g(u)\partial_tu\in L^{p}(0,T;(W^{1,q}(\Omega))')$ for $1\leq p<2$ and $q>2$.
\item[]
\item[iii)] $u(x,0)=u_0(x)$ for all $x\in \Omega$,
\end{itemize}
which solves \eqref{eqn:u}-\eqref{eqn:mu} in the following weak sense
\begin{itemize}
\item[a)]There exists a set $B\subset \Omega_T$ with $|\Omega_T\backslash B|=0$ and a function $\zeta:\Omega_T\rightarrow \R^n$ satisfying $\chi_{B\cap P}M(u)\zeta\in L^{\frac{p}{p-1}}(0,T;L^{\frac{q}{q-1}}(\Omega,\R^n))$ such that
          \begin{equation}
           \int_0^T<g(u)\partial_tu,\phi>_{(W^{1,q}(\Omega))',W^{1,q}(\Omega)}dt=-\int_{B\cap P}M(u)\zeta \cdot \nabla \phi dxdt \label{eqn:dgrrv}
          \end{equation} for all $\phi \in L^p(0,T;W^{1,q}(\Omega))$ with $p,q>2$. Here $P:=\{(x,t)\in \Omega_T: |1-u^2|\neq 0\}$ is the set where $M(u), g(u)$ are nondegenerate and $\chi_{B\cap P}$ is the characteristic function of set $B\cap P$.

\item[b)] Assume $u\in L^2(0,T;H^2(\Omega)).$  For any open set $U\in \Omega_T$ on  which $g(u)>0 $ and $\nabla \Delta u \in L^p(U)$ for some $p>1$, we have
\begin{equation}
\zeta=\frac{-\nabla \Delta u+q''(u)\nabla u+\nabla (-\Delta)^{\frac{1}{2}}u}{g(u)}-\frac{g'(u)}{g^2(u)}\left(-\Delta u+q'(u)+(-\Delta)^{\frac{1}{2}} u \right)\nabla u . \label{eqn:zetaUrrv}
\end{equation} a.e  in $U$.
\end{itemize}
Moreover, the following energy inequality holds for all $t>0$
\begin{eqnarray}
\notag &&\int_{\Omega}\left(\frac{1}{2}|\nabla u(x,t)|^2+q(u(x,t))\right)dz+\int_{\Omega_r\cap B\cap P}M(u(x,\tau))|\zeta(x,\tau)|^2 dxd\tau\\
\label{eqn:ineqrrv}&\leq &\int_{\Omega}\left(\frac{1}{2}|\nabla u_0(x)|^2+q(u_0(x))\right)dx.
\end{eqnarray}
\end{theorem}

\textbf{ACKNOWLEDGEMENTS}
 X.H. Niu's research is supported by National Natural Science Foundation of China under the grant number 11801214 and the Natural Science Foundation of Fujian Province of China under the grant number 2021J011193.
 Y. Xiang's research is supported by the Hong Kong Research Grants Council General Research
Fund 16307319.
 X. Yan's research is supported by a Research Excellence Grant and CLAS Dean's Summer Research Grant from University of Connecticut.

\FloatBarrier
\bibliographystyle{siam}
\bibliography{references}

\end{document}